\documentstyle{article}
\newtheorem{theorem}{Theorem}

\newtheorem{definition}[theorem]{Definition}
\newtheorem{example}[theorem]{Example}
\newtheorem{remark}[theorem]{Remark}
\newtheorem{corollary}[theorem]{Corollary}

\newcommand{\examp}[1]
  {\begin{example} {\rm #1} \end{example}}
\newcommand{\rem}[1]
  {\begin{remark} {\rm #1} \end{remark}}

\def\QED{\quad\blackslug\lower 8.5pt\null}

\newcommand{\crazy}[2]{\displaystyle{\mathop{#1}_{#2}}
\vphantom{\displaystyle{#1}}}

\begin{document}

\begin{center}
{\Large\bf A CLASSIFICATION AND EXAMPLES}

 \vspace*{2mm}

{\Large\bf OF FOUR-DIMENSIONAL ISOCLINIC  }

\vspace*{2mm}

{\Large\bf   THREE-WEBS}

\vspace*{4mm}

{\large Vladislav V. Goldberg}
\end{center}

 Abstract. {\footnotesize  A classification and examples  of four-dimensional
 isoclinic three-webs
of codimension two are given. The examples considered prove the existence
theorem for many classes of webs for which the general existence
theorems are not proved yet.}

 \vspace*{8mm}

\setcounter{section}{-1}

\section{Introduction} In the 1980s while studying rank problems
for webs (see the papers [G 83, 92, 93], the review paper [AG 99a]
and the monograph [G 88], Ch. 8), the author has constructed three
examples of exceptional
 four-webs $W (4, 2, 2)$ of maximum 2-rank on a
four-dimensional manifold  $X^4$ (see [G 85, 86, 87], [AG 99a], and
the books [G 88], [AS 92], and [AG 96]). They are exceptional
since  they are of maximum rank but not  algebraizable.
 H\'{e}naut [H 98]  named these webs after the author
and denoted them by $G_1 (4, 2, 2), G_2 (4, 2, 2)$, and $G_3 (4, 2, 2)$.

When the author was constructing these examples of four-webs, he
 considered numerous examples of three-webs $W (3, 2, 2)$ on $X^4$ and
proved that  three of them can be  expanded  to exceptional
 four-webs $W (4, 2, 2)$ of maximum 2-rank (see [G 85, 86, 87, 88]).

  However, it appeared that many  examples
of three-webs $W (3, 2, 2)$ that the author has constructed in that research
 are useful when one studies different classes of multidimensional
three-webs (isoclinic, hexagonal, transversally
geodesic, algebraizable, Bol's webs, etc.). Some of these examples
of  webs $W (3, 2, 2)$ were published in
the author paper [G 92] and the book [G 88], and some were
included in problem sections of the book [AS 92].

The aim of the current paper is to present some of these examples both
published and unpublished following some classification for them.
 We will include into this paper some of three-webs $W (3, 2, 2)$  considered by
Bol [B 35] and Chern [C 36]. These examples were already
considered in   [G 85, 86, 87, 88] from the rank point of view. It
is worth to consider them again in order to provide an up-to-date
characterization of webs by indicating to which classes
 a web belongs (isoclinic, transversally geodesic,
hexagonal, algebraizable, group, parallelizable, etc.). Most of
these classes were not known in the 1930s.

The examples mentioned above prove  existence
 of many classes of webs for which the general existence theorem
is not proved yet.

In this paper we will consider only isoclinic webs $W (3, 2, 2)$.
A classification and examples of nonisoclinic webs $W (3, 2, 2)$
will be presented elsewhere.

 \section{The transversal distribution of a web {\protect\boldmath\(W (3, 2, 2)\)
\protect\unboldmath}}

{\bf 1.} The leaves of the foliation $\lambda_\xi, \; \xi =
1, 2, 3$,  of a web $W (3, 2, 2)$ are determined by the equations
$\crazy{\omega}{\xi}^i = 0, \; i = 1, 2$,   where
\begin{equation}\label{eq:1}
\crazy{\omega}{1}^i + \crazy{\omega}{2}^i + \crazy{\omega}{3}^i = 0
\end{equation}
(see, for example, [G 88], Section {\bf 8.1} or
[AS 92], Section  {\bf 1.3}).

The structure equations of such a web can be written in the form
\begin{equation}\label{eq:2}
\renewcommand{\arraystretch}{1.3}
 \left\{
\begin{array}{ll}
   d \crazy{\omega}{1}^i =
\crazy{\omega}{1}^j \wedge \omega_j^i +  a_j \crazy{\omega}{1}^j
\wedge \crazy{\omega}{1}^i, \\  d\crazy{\omega}{2}^i =
\crazy{\omega}{2}^j \wedge \omega_j^j - a_j \crazy{\omega}{2}^j
\wedge \crazy{\omega}{2}^i.
\end{array}
      \right.
\renewcommand{\arraystretch}{1}
 \end{equation}
 The differential prolongations of equations (2) are
(see [G 88], Sections {\bf 8.1} and {\bf 8.4} or [AS 92], Section
{\bf 3.2}):
\begin{equation}\label{eq:3}
  d\omega_j^i - \omega_j^k \wedge
\omega_k^i =   b_{jkl}^i \crazy{\omega}{1}^k \wedge
\crazy{\omega}{2}^l,
 \end{equation}
\begin{equation}\label{eq:4}
da_i - a_j \omega_i^j = p_{ij} \crazy{\omega}{1}^j  +  q_{ij}
\crazy{\omega}{2}^j,
\end{equation}
 where
 \begin{equation}\label{eq:5}
b^i_{[j|l|k]} = \delta^i_{[k} p_{j]l}, \;\; b^i_{[jk]l} =
\delta^i_{[k} q_{j]l}.
\end{equation}
The quantities
\begin{equation}\label{eq:6}
 a_{jk}^i = a_{[j} \delta_{k]}^i
 \end{equation}
 and $b^i_{jkl}$ are the {\em
torsion and curvature tensors} of a three-web \linebreak $W (3,
2, 2)$. Note that for webs $W (3, 2, 2)$ the torsion tensor
$a^i_{jk}$ always has structure (6), where  $a = \{a_1, a_2\}$ is
a covector. If $a = 0$, then a web $W (3, 2, 2)$ is {\em isoclinicly
geodesic}. In what follows in this section, {\em we will assume that $a
\neq 0$, i.e., a web $W (3, 2, 2)$ is nonisoclinicly geodesic.}

{\bf 2}. For a web $W (3, 2, 2)$, a transversally geodesic
distribution is defined (cf. [AS 92], Section {\bf 3.1}) by the
equations
$$ \xi^2 \crazy{\omega}{1}^1 - \xi^1
\crazy{\omega}{1}^2= 0, \;\; \xi^2 \crazy{\omega}{2}^1 - \xi^1
\crazy{\omega}{2}^2 = 0.
$$
If we take $\displaystyle
\frac{\xi^1}{\xi^2} = - \frac{a_2}{a_1}$, we obtain an
invariant transversal distribution $\Delta$ defined by the
equations
\begin{equation}\label{eq:7}
 a_1 \crazy{\omega}{1}^1 + a_2 \crazy{\omega}{1}^2 = 0, \;\; a_1
 \crazy{\omega}{2}^1 + a_2 \crazy{\omega}{2}^2 = 0.
  \end{equation}

We will call the distribution $\Delta$ defined by equat\-ions (7)
the {\em transversal $a$-distribution} of a web $W (3, 2, 2)$
since it is defined by the covector $a$.
Note that for isoclinicly geodesic webs $W (3, 2, 2)$, for which
$a_1 = a_2 =0$, the transversal distribution is not defined.

The following theorem gives the conditions of integrability of the
distribution $\Delta$ (see the proofs of this and other results
of this section in the forthcoming paper [AG 99b]).

\begin{theorem}
The transversal $a$-distribution $\Delta$ defined by
 equa\-tions $(7)$  is integrable if and only if the components
$a_1$ and $a_2$ of the covector $a$ and their Pfaffian derivatives
$p_{ij}$ and $q_{ij}$ satisfy the conditions
\begin{equation}\label{eq:8}
\renewcommand{\arraystretch}{1.3}
\left\{
\begin{array}{ll}
a_2^2 p_{11} -  2a_1 a_2 p_{(12)} + a_1^2 p_{22} = 0, \\
a_2^2 q_{11} -  2a_1 a_2 q_{(12)} + a_1^2 q_{22} = 0.
        \end{array}
\right.
\renewcommand{\arraystretch}{1}
\end{equation}
\end{theorem}

Note that for a web $W (3, 2, 2)$, it is always possible
to take a specialized frame in which  there is a relation
between the components $a_1$ and $a_2$  of the covector $a$.
For example, if the transversal distribution $\Delta$ coincides
with the distribution $\crazy{\omega}{\alpha}^1 = 0$ or $\crazy{\omega}{\alpha}^2 = 0$
or $\crazy{\omega}{\alpha}^1 + \crazy{\omega}{\alpha}^2 = 0$,
 then we have $a_2 = 0$ or $a_1 = 0$ or $a_1 = a_2$, respectively.
Note that in these cases
the forms $\omega_2^1, \; \omega_1^2$, and  $\omega_1^1 + \omega_1^2
- \omega_2^1 - \omega_2^2$,
respectively, are expressed in terms of
the basis forms $\crazy{\omega}{\alpha}^i, \; \alpha = 1, 2$,
i.e., in these cases we have
$$
\pi_2^1 = 0, \;\; \pi_1^2 = 0, \;\; \pi_1^1 + \pi_1^2 - \pi_2^1 - \pi_2^2 = 0,
$$
respectively, where $\pi_i^j = \omega_i^j\Bigl|_{\crazy{\omega}{\alpha}^i = 0}$.

 In examples, that we are going to present in this paper, such situations
 will occur. So, we present here the conditions of integrability for these
 three cases.

\begin{corollary}
If for a web $W (3, 2, 2)$ one of the following conditions
\begin{equation}\label{eq:9}
a_2 = 0, \;\; \pi_2^1 = 0,
\end{equation}
\begin{equation}\label{eq:10}
a_1 = 0, \;\; \pi_1^2 = 0,
\end{equation}
\begin{equation}\label{eq:11}
a_1 = a_2, \;\; \pi_1^1 + \pi_1^2 - \pi_2^1 - \pi_2^2 = 0
\end{equation}
holds, then the $a$-distribution $\Delta$ coincides with
the distribution $\crazy{\omega}{\alpha}^1 = 0, \;  \crazy{\omega}{\alpha}^2 = 0$
or   $\crazy{\omega}{\alpha}^1 + \crazy{\omega}{\alpha}^2 = 0$, respectively.
This $a$-distribution is integrable if and only if the quantities
$p_{ij}$ and $q_{ij}$ satisfy  respectively the following
conditions:
\begin{equation}\label{eq:12}
p_{22} = q_{22} = 0,
\end{equation}
\begin{equation}\label{eq:13}
p_{11} = q_{11} = 0,
\end{equation}
\begin{equation}\label{eq:14}
p_{11} - 2 p_{(12)} + p_{22} = 0, \;\;
q_{11} - 2 q_{(12)} + q_{22} = 0.
\end{equation}
\end{corollary}

Each of relations (8), (12), (13), and (14) gives two
conditions which Pfaffian derivatives $p_{ij}$ and $q_{ij}$ of the co-vector $a$ must
satisfy in order for the $a$-distribution $\Delta$
of a web $W (3, 2, 2)$ to be integrable.

{\bf 3.} On surfaces $V^2$ of the integrable
$a$-distribution $\Delta$, foliations $\lambda_\xi$ of a web
$W (3, 2, 2)$ cut two-dimensional three-webs $W (3, 2, 1)$.
It is well-known (see, for example, [AS 92], Section {\bf 3.1}),
that if a web $W (3, 2, 2)$ is
hexagonal, then it is transversally geodesic, and
 two-dimensional  three-webs on all integral surfaces $V^2$ of the
transversally geodesic distribution are hexagonal, and conversely.
However, it is not true for the integrable transversal distribution
$\Delta$: the condition of hexagonality of all webs $W (3, 2, 1)$
is weaker than the condition $b^i_{(jkl)} = 0$ of hexagonality of
the web $W (3, 2, 2)$.

 The following theorem describes
   webs $W (3, 2, 2)$ for which hexagonality of all
  two-dimensional \mbox{webs}  $W (3, 2, 1)$  cut by the foliations of webs
  $W (3, 2, 2)$ on the integral surfaces $V^2$ of their  transversal distribution
$\Delta$ implies  hexagonality of a web $W (3, 2, 2)$.

\begin{theorem}
\begin{description}
\item[(i)] Let $W (3, 2, 2)$ be a web with nonvanishing covector  $a \neq 0$ and
with integrable transversal $a$-distribution $($conditions $(8)$ hold$)$.
The lines of all two-dimensional  webs $W (3, 2, 1)$ cut by the foliations of
a web $W (3, 2, 2)$  with constant $t$ on the integral surfaces $V^2$
of its transversal $a$-distribution $\Delta$ are geodesic in all
affine connections defined by the connection forms $\theta_j^i =
\omega_j^i + a^i_{jk} (p \crazy{\omega}{1}^k + q
\crazy{\omega}{2}^k)$ if and only if the equation
\begin{equation}\label{eq:15}
\omega_2^1 = t^2 \omega_1^2 + t (\omega_1^1 - \omega_2^2)
\end{equation}
holds, where $t = \displaystyle \frac{a_2}{a_1}$ is an arbitrary
constant.

 \item[(ii)] A web $W (3, 2, 2)$ satisfying equation $(15)$ is  hexagonal
if and only if for any constant $t$ all  two-dimensional  webs $W (3, 2, 1)$
cut by the foliations of  a
web $W (3, 2, 2)$  on the integral surfaces $V^2$
of its transversal $a$-distribution $\Delta$ are hexagonal.
The condition of integrability of $\Delta$ implies the equation
\begin{equation}\label{eq:16}
\renewcommand{\arraystretch}{1.3}
\begin{array}{ll}
b^2_{111} t^4 \!\!\!\!\!&- (3b^2_{(112)} - b^1_{111}) t^3 + 3(b^2_{(122)}
- b^1_{(112)}) t^2 \\
\!\!\!\!\!& - (b^2_{222} - 3b^1_{(122)}) t - b^1_{222} = 0,
  \end{array}
\renewcommand{\arraystretch}{1}
\end{equation}
and the condition of hexagonality of  webs  $W (3, 2, 1)$ is
\begin{equation}\label{eq:17}
- b^1_{111} t^3 + 3b^1_{(112)} t^2 - 3 b^1_{(122)} t + b^1_{222} = 0
\end{equation}
or
\begin{equation}\label{eq:18}
- b^2_{111} t^3 + 3b^2_{(112)} t^2 - 3 b^2_{(122)} t + b^2_{222} = 0.
\end{equation}
\end{description}
\end{theorem}

Note that conditions (15) and (8) imply that $t = \mbox{const.}$,
and any two of
conditions  $(16)$, $(17)$, and $(18)$ imply the third one. Conditions
$(17)$ and $(18)$ show that if they are valid for {\em any constant} $t$,
then $b^i_{(jkl)} = 0$, and a web $W (3, 2, 2)$ is hexagonal.

We consider now three cases that are similar to those in Corollary 2
with one exception: instead of having $\pi_2^1 = 0, \;\pi_1^2 = 0$, and
$\pi_1^1 + \pi_1^2 - \pi_2^1 - \pi_2^2 = 0$,
we request the stronger conditions $\omega_2^1 = 0, \; \omega_1^2 = 0$, and
$\omega_1^1 + \omega_1^2 - \omega_2^1 - \omega_2^2 = 0$, respectively.

\begin{corollary}
If for a web $W (3, 2, 2)$ one of the following conditions
\begin{equation}\label{eq:19}
  a_2 = 0, \; \omega_2^1 = 0,
\end{equation}
\begin{equation}\label{eq:20}
  a_1 = 0, \; \omega_1^2 = 0,
\end{equation}
\begin{equation}\label{21}
  a_1 = a_2, \;\; \omega_1^1 + \omega_1^2 = \omega_2^1 + \omega_2^2,
\end{equation}
holds,  then the  transversal $a$-distribution
$\Delta$ is the  distribution  defined by  the equations
 $\crazy{\omega}{\alpha}^1 = 0, \; \crazy{\omega}{\alpha}^2 = 0$,
 or  $\crazy{\omega}{\alpha}^1 + \crazy{\omega}{\alpha}^2 = 0$,
 respectively. This distribution  is integrable,   and such a web
is foliated into  two-dimensional  hexagonal webs $W (3, 2, 1)$ belonging to
 integral surfaces of $\Delta$ if and only if
\begin{equation}\label{eq:22}
 b^1_{222} = 0, \;\; b^2_{222} = 0,
\end{equation}
\begin{equation}\label{eq:23}
 b^1_{111} = 0, \;\; b^2_{111} = 0,
\end{equation}
\begin{equation}\label{eq:24}
\renewcommand{\arraystretch}{1.3}
\left\{
\begin{array}{ll}
- b^1_{111} + 3 (b^1_{(112)} - b^1_{(122)}) + b^1_{222} = 0,  \\
- b^2_{111} + 3 (b^2_{(112)} - b^2_{(122)}) + b^2_{222} = 0,
  \end{array}
\right.
\renewcommand{\arraystretch}{1}
 \end{equation}
 respectively.
\end{corollary}

\section{The isoclinic webs $W (3, 2, 2)$}

{\bf 1.}
A one-parameter family of isoclinic bivectors $\zeta_1 \wedge \zeta_2$, where
$\zeta_i = \zeta^1 \crazy{e}{2}_i -  \zeta^2 \crazy{e}{1}_i$,
is associated with a web $W (3, 2, 2)$ (see [AS 92], Sections {\bf 1.3}
and {\bf 3.2} or [G 88], Section {\bf 1.11}). A surface is {\em isoclinic} if
it is tangent to an isoclinic bivector at each of its points.
A web $W (3, 2, 2)$ is called {\em isoclinic} if a one-parameter family
of isoclinic surfaces passes through any point of this web.

It is proved (see [G 88], Section {\bf 8.1}, [G 92] or [AS 92].
Section {\bf 3.2}) that {\em  a web $W (3, 2, 2)$ is
isoclinic if and only if the quantities
$p_{ij}$ and $q_{ij}$ in equations $(4)$ are symmetric}:
\begin{equation}\label{eq:25}
p_{ij} = p_{ji}, \;\; q_{ij} = q_{ji}.
\end{equation}
  Note that for webs $W (3, 2, r), \; r > 2$, the structure (6) of its
 torsion tensor is necessary and sufficient for its isoclinicity.
 As we noted earlier,  the torsion tensor of a web $W (3, 2, 2)$ always has form (6), and
  conditions (25) are necessary and sufficient for this web to be isoclinic.

For an isoclinic web $W (3, 2, 2)$, we can take
\begin{equation}\label{26}
  p_{ij} = f_{ij} - h_{ij}, \;\;   q_{ij} = g_{ij} - h_{ij},
\end{equation}
where $f_{ij}, g_{ij}, h_{ij}$ are symmetric (0, 2)-tensors.
Then the curvature tensor of the web $W (3, 2, 2)$ has the form
\begin{equation}\label{eq:27}
 b_{jkl}^i = a_{jkl}^i + f_{jk}  \delta_l^i +
 g_{lj}  \delta_k^i + h_{kl}  \delta_j^i,
\end{equation}
where
\begin{equation}\label{28}
  a^i_{jkl} = b^i_{(jkl)} - (f_{(jk} + g_{(jk} + h_{(jk}) \delta_{l)}^i,
\end{equation}
and
\begin{equation}\label{eq:29}
 a_{ikl}^i = 0
\end{equation}
(see [AS 92], Section {\bf 3.2}).
The condition (29) allows us to find the tensor $h_{ij}$:
\begin{equation}\label{eq:30}
 h_{ij}  = \frac{3}{4} b_{(kij)}^k - f_{ij} - g_{ij} .
\end{equation}
So, now the tensors  $f_{ij}, g_{ij}$, and $h_{ij}$ are
well-defined by equations (30) and (26).
It follows from (30) and (26) that the tensor $h_{ij}$ can be also
expressed as
\begin{equation}\label{eq:31}
h_{ij} = \frac{1}{4} b_{(kij)}^k - \frac{1}{3} (p_{ij} + q_{ij}).
\end{equation}

{\bf 2.} We will combine here some  results giving an analytic
characterizations of the most important classes of isoclinic webs
$W (3, 2, 2)$ (see [AS 92] and [G 88] for details). Note that
for the almost algebraizable, almost Bol, and almost parallelizable webs
these analytic characterizations have been used in [G 86, 87] (see also [G 88],
Ch. 8) as their definitions.

\vspace*{3mm}

\begin{theorem}
 For an isoclinic web $W (3, 2, 2)$
of different types we have the following analytic
characterizations:
\begin{description}
\item[(a)]
  An  isoclinic
web $W (3, 2, 2)$ is isoclinicly geodesic  if and only if
\begin{equation}\label{eq:32}
a_1 = a_2  = 0.
\end{equation}
\item[(b)]
  An  isoclinic web $W (3, 2, 2)$ is transversally
geodesic $($and therefore \newline Grassmannizable$)$ if and only if
\begin{equation}\label{eq:33}
 a_{jkl}^i = 0.
\end{equation}
\item[(c)]
  An  isoclinic web $W (3, 2, 2)$ is almost algebraizable
 if and only if
\begin{equation}\label{eq:34}
 f_{ij} + g_{ij} + h_{ij} = 0.
\end{equation}
\item[(d)]
  An  isoclinic web $W (3, 2, 2)$ is almost Bol web $B_m$
 if and only if
\begin{equation}\label{eq:35}
 f_{ij} + g_{ij} = 0, \;\; h_{ij} = 0.
\end{equation}
\item[(e)]
  An  isoclinic web $W (3, 2, 2)$ is almost parallelizable
 if and only if
\begin{equation}\label{eq:36}
 f_{ij} = g_{ij} = h_{ij} = 0.
\end{equation}
\item[(f)]
   An isoclinic web $W (3, 2, 2)$ is hexagonal
 $($and therefore algebraizable$)$ if and only if
\begin{equation}\label{eq:37}
a^i_{jkl} = 0, \;\; f_{ij} + g_{ij} + h_{ij} = 0.
\end{equation}
\item[(g)]
   An isoclinic web $W (3, 2, 2)$ is a
 Bol  web  $B_m$  if and only if
\begin{equation}\label{eq:38}
a^i_{jkl} = 0, \;\; f_{ij} + g_{ij} = 0, \;\; h_{ij} = 0.
\end{equation}
%\item[e)]
% An isoclinic web $W (3, 2, 2)$ is a group
%web if and only if
%\begin{equation}\label{eq:?}
% a_{jkl}^i + f_{jk}  \delta_l^i +
% g_{lj}  \delta_k^i + h_{kl}  \delta_j^i = 0,
%% a_{jkl}^i +  \delta_{(j}^i (f_{kl) + g_{kl} + h_{kl}) = 0.
%\end{equation}
\item[(h)]
  An isoclinic web $W (3, 2, 2)$ is a   group web if and
only if
\begin{equation}\label{eq:39}
a^i_{jkl} = 0, \;\; f_{ij} = g_{ij} = h_{ij} = 0.
\end{equation}
\item[(i)]
  An isoclinic   web $W (3, 2, 2)$ is parallelizable if and
only if conditions $(32)$ and $(39)$ hold.
\end{description}
\end{theorem}

{\sf Proof}. First note that the term "almost" is used
in parts (c), (d), and (e) of the theorem since a web $W (3, 2, 2)$
is not assumed to be transversally geodesic.

 We will prove only parts (h) and (i). The proof of
parts (a)--(g)  can be found in  [A 74]; [AS 92], Section {\bf 3.3};
[G 92] and [G 88], Section {\bf 8.1}.

{\em  Part {\rm (h)}}. Since a group web is transversally
geodesic, we have (33). Since it is  a group web, we also have
$b^i_{jkl} = 0$, and this and (27) imply that
$$
f_{jk}  \delta_l^i +  g_{lj}  \delta_k^i + h_{kl}  \delta_j^i = 0.
$$
Contracting the
last equations with respect to the indices $i$ and $j$, $i$ and
$k$, $i$ and $l$, we find that
$$
2f_{ij} + g_{ij} + h_{ij} = 0, \;\; f_{ij} + 2 g_{ij} + h_{ij} = 0, \;\;
f_{ij} + g_{ij} + 2h_{ij} = 0.
$$
These equations imply (39).
Thus, equations (39) are necessary and sufficient for an isoclinic
web to be a  group web.

{\em   Part {\rm (i)}}. Since a parallelizable web is a group web,
we have equations (39). In addition, the torsion tensor $a^i_{jk}$
of such a web must vanish, i.e., the web is isoclinicly
geodesic.  By (6), this implies (32).
 \rule{3mm}{3mm}

{\bf 3.} Suppose that three foliations $\lambda_1 ,\; \lambda_2$,
and $\lambda_3$ of a web $W (3, 2, 2)$  are given as the level
sets $u_{\xi}^i = \mbox{const.}\; (\xi = 1, 2, 3$), of the
following functions:
\begin{equation}\label{eq:40}
 \lambda_1 : u_1^i = x^i ;\;\;
 \lambda_2 : u_2^i = y^i ;\;\;
 \lambda_3 : u_3^i = f^i(x^j,y^k),\;\;\;\;
 i,j,k = 1,2 .
\end{equation}

In order to characterize a web $W (3, 2, 2)$ given by
(40), we must find the forms $\crazy{\omega}{\alpha}^i, \;\alpha = 1, 2, \;
\omega_j^i$, and the functions $a_{jk}^i ,\; b_{jkl}^i,\; a_i ,\;
a_{jkl}^i ,\; f_{ij}, g_{ij}, h_{ij}$.

The forms $\crazy{\omega}{\alpha}^i ,\; \omega_j^i$ and the
functions $a_{jk}^i$ and $b_{jkl}^i$ can be found my means of the
following formulae (see [AS 71], or [G 88], Section {\bf 8.1}, or
[AS 92], Section {\bf 1.6}):
\begin{equation}\label{eq:41}
 \crazy{\omega}{1}^i = \bar{f}_j^i dx^j ,\;\;\;
 \crazy{\omega}{2}^i = \tilde{f}_j^i dy^j ,\;\;\;
 \crazy{\omega}{3}^i = - du_3^i ,
\end{equation}
where $$
 \bar{f}_j^i = \displaystyle
\frac{\partial{f^i}}{\partial{x^j}}, \;\;\;\;
 \tilde{f}_j^i = \displaystyle
\frac{\partial{f^i}}{\partial{y^j}} , \;\;\;\;
 \det{(\bar{f}_j^i)} \neq 0 , \;\;\;
 \det{ (\tilde{f}_j^i)} \neq 0,
$$ and
\begin{equation}\label{eq:42}
 d\crazy{\omega}{1}^i = - d\crazy{\omega}{2}^i = \Gamma_{jk}^i
 \crazy{\omega}{1}^j \wedge \crazy{\omega}{2}^k ,
\end{equation}
\begin{equation}\label{eq:43}
 \Gamma_{jk}^i = - \displaystyle
\frac{\partial^2 f^i}{\partial{x^l} \partial{y^m}}
 \bar{g}_j^l \tilde{g}_k^m ,
\end{equation}
\begin{equation}\label{eq:44}
 \omega_j^i = \Gamma_{kj}^i \crazy{\omega}{1}^k + \Gamma_{jk}^i
 \crazy{\omega}{2}^k ,
\end{equation}
\begin{equation}\label{eq:45}
 a_{jk}^i = \Gamma_{[jk]}^i ,
\end{equation}
\begin{eqnarray}\label{eq:46}
\renewcommand{\arraystretch}{1.3}
 b_{jkl}^i & =  \displaystyle \frac{1}{2} \Biggl(
 \displaystyle
\frac{\partial{\Gamma_{kl}^i}}{\partial{x^m}} \bar{g}_j^m +
\displaystyle
 \frac{\partial{\Gamma_{jl}^i}}{\partial{x^m}} \bar{g}_k^m
- \displaystyle
 \frac{\partial{\Gamma_{kj}^i}}{\partial{y^m}} \tilde{g}_l^m
- \displaystyle \frac{\partial{\Gamma_{kl}^i}}{\partial{y^m}}
\tilde{g}_j^m \nonumber \\ &+   \Gamma_{jl}^m \Gamma_{km}^i -
\Gamma_{kj}^m \Gamma_{ml}^i
 + 2\Gamma_{kl}^m a_{mj}^i\Biggr) .
\renewcommand{\arraystretch}{1}
\end{eqnarray}

As to the functions $a_i ,\; f_{ij}, g_{ij}, h_{ij}$, and
$a_{jkl}^i$, they can be easily calculated from equations (2),
(4), (6), (26), and (30).

In what follows, we will always assume that three  foliations of a
web $W (3, 2, 2)$ are given by
\begin{equation}\label{eq:47}
\renewcommand{\arraystretch}{1.3}
\left\{
        \begin{array}{ll}
         \lambda_1 : x^1 = \mbox{const.},\;\;\; x^2 =
\mbox{const.};\;\;\; \\
         \lambda_2 : y^1 = \mbox{const.},\;\;\; y^2 =
\mbox{const.}; \\ u^1_3 = f^1 (x^j, y^k) = \mbox{const.},\;\;\;
u_3^2 = f^2 (x^j, y^k) =\mbox{const.}
        \end{array}
        \right.
\renewcommand{\arraystretch}{1.3}
\end{equation}
In examples of webs, that we are going to present
in Section {\bf 3},  we will  only specify the functions $f^1
(x^j, y^k)$ and $f^2 (x^j, y^k)$.

Note that it follows from (44) that all forms $\omega_j^i$ are
expressed in terms of the forms $\crazy{\omega}{1}^i$ and
$\crazy{\omega}{2}^i$ only. This means that the forms $\pi_j^i =
\omega_j^i\Bigl|_{\crazy{\omega}{\alpha}^i = 0}$
 vanish, $\pi_j^i = 0$.
  It follows that the second conditions in equations
  (9), (10), and  (11) are always valid for webs
  defined by equations (47), and the meaning of
  equations  (9), (10), and (11) is that the transversal
  distribution $\Delta$ coincides with the distribution
  $\crazy{\omega}{1}^i = 0$ or $\crazy{\omega}{2}^i = 0$, or
   $\crazy{\omega}{1}^i  + \crazy{\omega}{2}^i = 0$, respectively.
%   or $\crazy{\omega}{1}^i - \crazy{\omega}{2}^i = 0$.
 %all the  tensors associated with a web $W (3, 2, 2)$ defined
%by equations (60) are absolute tensors, and each of their
%component is an absolute invariant. Thus the vanishing of any of
%these components distinguishes a class of webs $W (3, 2, 2)$
%defined by equations (58).

{\bf 4.} We will present a classification of isoclinic webs $W (3, 2, 2)$
given by equations (47).
\begin{description}
\item[{\bf A.}] Nonisoclinicly geodesic webs ($a_1^2 +  a_2^2 > 0$).
  \begin{description}  \item[{\bf A${}_{1}.$}] Webs  with the
             integrable  transversal distribution $\Delta$ ((8) holds).
     \begin{description} \item[{\bf A${}_{11}.$}] Webs  with the
             integrable  transversal distribution $a_1 \crazy{\omega}{1}^i
+ a_2 \crazy{\omega}{2}^i = 0, \; a_1, a_2 \neq 0, \;
 t = \frac{a_2}{a_1} = \mbox{{\rm const.}}$
                  ((8) holds).
        \begin{description}
                                          \item[{\bf A${}_{111}.$}]
                                          Webs   foliated into planar
                                          hexagonal webs $W (3, 2, 1)$
                                             ((17) and (18) hold).
                                              \item[{\bf A${}_{112}$.}]
                                           Webs with  the integrable
                                           transversal distribution
    $\crazy{\omega}{\alpha}^1 + \crazy{\omega}{\alpha}^2   = 0$ ((11) and (14) hold).
         \begin{description}
                                                     \item[{\bf A$_{1121}$.}]
                                                     Webs   foliated into
                                                     planar hexagonal webs
                                                     $W (3, 2, 1)$ ((21) and (24) hold).
         \end{description}
\end{description}

%    \item[{\bf A${}_{113}$.}]
% Webs with the integrable  transversal distribution
% $\crazy{\omega}{1}^i  - \crazy{\omega}{2}^i = 0$
% ((19) and (20) hold).
%\begin{description}
%                    \item[{\bf A${}_{1131}$.}] Webs   foliated
%into planar hexagonal webs $W (3, 2, 1)$ ((43) and (44) hold).
        %\end{description}

           \item[{\bf A${}_{12}.$}] Webs  with the
            integrable  transversal distribution $\crazy{\omega}{\alpha}^1 = 0$
                  ((9) and (12) hold).
  \begin{description}
                    \item[{\bf A${}_{121}.$}] Webs   foliated
into planar hexagonal webs $W (3, 2, 1)$
 ((19) and (22) hold).
     \end{description}
            \item[{\bf A${}_{13}$.}] Webs  with the integrable  transversal
                                     distribution $\crazy{\omega}{\alpha}^2 = 0$
                                     ((10) and (13) hold).
   % \end{description}

   \begin{description}
                                       \item[{\bf A${}_{131}.$}] Webs   foliated
                                        into planar  hexagonal webs $W (3, 2, 1)$
                                        ((20) and (23) hold).
  \end{description}
  %\end{description}\end{description}
  \end{description}
      \item[{\bf A${}_2.$}] Webs with nonintegrable  transversal
                          distribution $\Delta$.
    \end{description}
 %\end{description}
%\end{description}

\item[{\bf B.}] Isoclinicly geodesic webs ((32) holds).
\item[{\bf C.}] Nontransversally geodesic webs ((33) does not hold).
\begin{description}
\item[{\bf C${}_1$.}] Almost algebraizable webs ((34) holds).
\begin{description}
\item[{\bf C${}_{11}$.}] Almost Bol webs ((35) holds).
\item[{\bf C${}_{12}$.}] Almost parallelizable webs ((36) holds).
\end{description}
\item[{\bf C${}_2$.}] Webs that are not almost algebraizable
((34) does not hold).
\end{description}
\item[{\bf D.}] Grassmannizable webs ((33) holds).
\begin{description}
\item[{\bf D${}_1$.}] Nonhexagonal webs ((37) does not hold).
\item[{\bf D${}_2$.}] Hexagonal webs ((37) holds).
\begin{description}
\item[{\bf D${}_{21}$.}] Bol webs ((38)  holds).
\item[{\bf D${}_{22}$.}] Nongroup webs ((39) does not hold).
\item[{\bf D${}_{23}$.}] Group webs ((39) holds).
\begin{description}
\item[{\bf D${}_{231}$.}] Nonparallelizable webs
                        ((32) does not hold).
\item[{\bf D${}_{232}$.}] Parallelizable webs ((32)
 holds).
\end{description}
\end{description}
\end{description}
\item[{\bf E.}] Webs with different relations among $p_{ij}$ and
                   $q_{ij}$. We indicate the most common classes
                   of such webs.
\begin{description}
\item[{\bf E${}_1$.}] Webs with $p_{ij} = q_{ij} = 0$.
\item[{\bf E${}_2$.}] Webs with $p_{1i} = q_{ij} = 0; \;\;
 p_{22} \neq 0.$
\item[{\bf E${}_3$.}] Webs with $p_{ij} = q_{1i} = 0; \;\;
 q_{22} \neq 0$.
\item[{\bf E${}_4$.}] Webs with $p_{1i} = q_{1i} = 0; \;\;
 p_{22} \neq 0, \;\; q_{22} \neq 0$.
\begin{description}
\item[{\bf E${}_{41}$.}] Webs with $p_{22} = - q_{22} \neq 0$.
\end{description}
\item[{\bf E${}_5$.}] Webs with $p_{2i} = q_{ij} = 0; \;\;
 p_{11} \neq 0.$
\item[{\bf E${}_6$.}] Webs with $p_{ij} = q_{2i} = 0; \;\;
 q_{11} \neq 0$.
\item[{\bf E${}_7$.}] Webs with $p_{2i} = q_{2i} = 0; \;\;
 p_{11} \neq 0, \;\; q_{11} \neq 0$.
\begin{description}
\item[{\bf E${}_{71}$.}] Webs with $p_{11} = - q_{11} \neq 0$.
\end{description}
\item[{\bf E${}_8$.}] Webs with $p_{ij} + q_{ij} = 0$.
\end{description}

\item[{\bf F.}] Webs extendable to a web $W (4, 2, 2)$ of maximum 2-rank.
\begin{description}
\item[{\bf F${}_1$.}] Webs extendable to an
algebraizable web $W (4, 2, 2)$ of maximum 2-rank.
\item[{\bf F${}_2$.}] Webs extendable to an exceptional
(nonalgebraizable) web $W (4, 2, 2)$ of maximum 2-rank.
\end{description}
\item[{\bf G.}] Webs nonextendable to a web $W (4, 2, 2)$ of
 maximum 2-rank.
\end{description}

\rem{It is easy to give  characterizations of  webs of the classes
{\bf  E${}_{1}$}--{\bf  E${}_{8}$}:

\begin{itemize}

\item Class {\bf  E${}_{1}$}:  the covector $a$ is covariantly
      constant on  an entire web of this class.

\item Class {\bf  E${}_{2}$}: for webs of this class, the component $a_1$ is covariantly
      constant on  an entire web, and $a_2$ is covariantly
      constant on the three-dimensional
      distribution defined by the equation $\crazy{\omega}{1}^2 = 0$.

\item Class {\bf  E${}_{4}$}: for webs of this class, the component $a_1$ is covariantly
      constant on  an entire web, and the component $a_2$ is covariantly
      constant on the two-dimensional distribution defined by the equations
      $\crazy{\omega}{1}^2 = 0, \;  \crazy{\omega}{2}^2 = 0$.

\item Class {\bf  E${}_{41}$}: for webs of this class, the component $a_1$ is covariantly
      constant on  an entire web, and the component $a_2$ is covariantly
      constant on the three-dimensional distribution defined by the equations
      $\crazy{\omega}{1}^2 - \crazy{\omega}{2}^2 = 0$.

\item Class {\bf  E${}_{8}$}: for webs of this class,   the covector $a$ is covariantly
      constant on  the two-dimensional distribution defined by the equations
      $\crazy{\omega}{1}^j  - \crazy{\omega}{2}^j = 0, \;\; j = 1, 2$.

The characterization of webs of the classes {\bf  E${}_{3}$}
(and {\bf  E${}_{5}$}, {\bf  E${}_{6}$}), {\bf  E${}_{4}$},
{\bf  E${}_{41}$} is similar to that of the classes
{\bf  E${}_{2}$}, {\bf  E${}_{7}$}, {\bf  E${}_{71}$}, respectively.
\end{itemize}
}

\rem{
If a web is given by equation (1), then
the webs of the classes {\bf A${}_{11}$}, {\bf A${}_{112}$},
%{\bf A${}_{113}$},
{\bf A${}_{12}, $} {\bf A${}_{13}$} as well as webs of the classes
{\bf A${}_{111}$}, {\bf A${}_{1121}$},
%{\bf A${}_{1131}$},
{\bf A${}_{121}$}, {\bf A${}_{131}$} could be equivalent one to another.
They are different by a location of the integrable transversal distribution
$\Delta$. In fact, if there is no additional conditions on webs, then
by transformations $x^i = \phi^i (x^j), \; y^i = \psi^i (y^j)$, we
can make the integrable transversal distributions $\Delta$ of any two
of them to coincide. However, in our examples in Section {\bf 3},
we always have $\pi_i^j = 0$, and  above mentioned
specializations are impossible.
%---webs are given with this kind or other kind of
%specialization done.
In addition, in our examples,
there will be additional conditions on webs, and in general, under the above mentioned
transformation (if it would be possible), these additional conditions for the first transformed web
 do not coincide with the additional conditions for the second
 web. These are the reasons that in our classification, we consider
all above mentioned classes as different ones.
}
\rem{The classification presented above has some overlapping classes:
for example, the class {\bf A${}_1$} is a subclass of the class
{\bf E${}_1$}, and the class {\bf  D${}_{232}$} is a subclass of the class
{\bf B}.
}

\rem{In general, we must prove the existence theorem for all
the classes listed above. Such a theorem can be proved
for a web of general kind of each class using the well-known
Cartan's test or it can be proved by finding examples
of webs of these classes. The examples
of webs in Section {\bf 3}  prove existence of webs of most of the classes
of our classification.
}

\section{Examples of isoclinic  webs $W (3, 2, 2)$}

\setcounter{theorem}{0}

\examp{\label{examp:1} When we constructed an exceptional
 web $G_1 (4, 2, 2)$ (see [G 85, 86] or
[G 88], Example {\bf 8.4.2}), we  investigated the web $W (3, 2, 2)$
introduced by Bol (see [B 35], p. 462). This web $W(3, 2, 2)$
is defined by
\begin{equation}\label{eq:48}
 u^1_3 = x^1 + y^1,\;\;\;
u_3^2 = (x^2 + y^2) (y^1 - x^1)
\end{equation}
in a domain of ${\bf R}^4$ where $x^1 \neq y^1$. Note that this
example of the web $G_1 (4, 2, 2)$  was also included into the books
[AS 92], Section {\bf 3.5}, Problem {\bf 6} and  Ch. {\bf 8},
Problem {\bf 4};  [AG 96], Example {\bf 5.5.1}; and the review
paper [AG 99a], Section {\bf 6.4}).

Using (41)--(46), (26), (28), and (31), for the web (48), we find that
\begin{equation}\label{eq:49}
\renewcommand{\arraystretch}{1.3}
\left\{
\begin{array}{ll}
\Gamma_{ij}^1 = \Gamma_{22}^2 = 0, \;\;\; \Gamma_{11}^2 = \displaystyle
\frac{2(x^2 + y^2)}{x^1 - y^1},\;\;\; \Gamma_{21}^2
= - \Gamma_{12}^2 = \displaystyle \frac{1}{x^1 -y^1};  \\
\omega_1^2 = \displaystyle \frac{1}{x^1 -y^1} \Bigl[(x^2 +y^2)
(2\crazy{\omega}{1}^1 + \crazy{\omega}{2}^1) + \crazy{\omega}{1}^2
- \crazy{\omega}{1}^1\Bigr], \\ \omega_2^2 = - \displaystyle
\frac{1}{x^1 -y^1} \Bigl(\crazy{\omega}{1}^1 -
\crazy{\omega}{2}^1\Bigr), \;\; \omega_i^1 = 0,\\
%(x^1 - y^1) - (dx^2 - dy^2 ),\;\; \omega_2^2
%= - d\ln (x^1 - y^1) ;\\
 a_1 = \displaystyle \frac{2}{(y^1 - x^1)} ,\;\; a_2 = 0,\;\;
p_{2i} = q_{2i} = 0,\;\;  p_{11} = -q_{11} =\displaystyle
\frac{2}{(x^1 - y^1)^2}, \\ b_{ijk}^1 = b_{222}^2 = b_{211}^2 =
b_{122}^2 = b_{212}^2 = b_{211}^2 = 0, \\ b_{112}^2 = - b_{121}^2
= \displaystyle \frac{2}{(x^1 - y^1)^2}, \;\;\;   b_{111}^2 =
-\displaystyle \frac{8}{(x^2 + y^2)(x^1 - y^1)^2},\\
          h_{ij} = 0, \;\; f_{11} =  -
          g_{11} =\displaystyle  \frac{2}{(x^1 - y^1)^2},\;\;\;
          f_{ij} = g_{ij} = 0 ,\;\;\; (i,j) \neq (1,1);\\
   a_{ijk}^1 = 0 ,\;\; a_{111}^2 = b_{111}^2 ,\;\; a_{ijk}^2 = 0,
          \;\;\; (i,j,k) \neq (1,1,1) .
         \end{array}
         \right.
\renewcommand{\arraystretch}{1}
\end{equation}

Equations (49) and (25)  show that {\em the three-web $(48)$
is isoclinic}. Since $a^2_{111} \neq 0$,   equations (49) and (33) show
that {\em this web is not transversally geodesic}. Since
we have
$$
 f_{ij} + g_{ij} = 0,\;\;\; h_{ij} = 0,
$$
the web  $(48)$ is an {\it almost Bol three-web} $(B_m)$ (see (35)). Thus, the
{\it it belongs to the class {\bf C${}_{11}$}}.  By (49), (19), and (22),
 {\it the web $(48)$ belongs to the class {\bf A$_{121}$}},
and {\em it also belongs to the class {\bf E$_{71}$}}.

{\em The web $(48)$ belongs to the class {\bf F$_2$}} since as was
 proved in [G 85, 86] (see also [G 88], Example {\bf 8.4.2}),
 {\em the three-web $(48)$ can be extended to an exceptional
web $G_1 (4, 2, 2)$ of maximum $2$-rank}. Note that
 there exists a one-parameter family of such extensions, and
the leaves of the fourth foliation of the web $W (4, 2, 2)$ are
defined as level sets of the functions
$$
\renewcommand{\arraystretch}{1.3}
\left\{
\begin{array}{ll}
 u_4^1 = \displaystyle
\frac{(x^1 - y^1)^2 (x^2 + y^2)^2}{(x^2 + C)(y^2 -C)},  \\ u_4^2 =
x^1 + y^1 + \displaystyle \frac{(y^1 - x^1)(x^2 +y^2)}{\sqrt{(x^2
+ C)(y^2 -C)}
 \cdot   \mbox{arctan} \;\; \sqrt{(y^2 - C)(x^2 + C)}}.
        \end{array}
\right.
\renewcommand{\arraystretch}{1}
$$
The only abelian equation of the exceptional webs $G_1 (4, 2, 2)$  of this
family has the form
$$
\Biggl(\displaystyle \frac{1}{x^2}\Biggr)
dx^1 \wedge dx^2 + \Biggl(\displaystyle \frac{1}{y^2} \Biggr) dy^1
\wedge dy^2 - \Biggl(\displaystyle\frac{1}{u_3^2}\Biggr)  du_3^1
\wedge du_3^2 - \Biggl(\displaystyle \frac{1}{2u_4^1}\Biggr)
du_4^1 \wedge du_4^2 = 0.
$$ }

\examp{\label{examp:2} When we constructed an exceptional
 web $G_2 (4, 2, 2)$ (see [G 86]; [G 88], Example {\bf 8.4.3}),
we  investigated the web $W (3, 2, 2)$ defined by the functions
\begin{equation}\label{eq:50}
 u_3^1 = x^1 + y^1, \;\;\;
         u_3^2 = -x^1 y^2 + x^2 y^1
\end{equation}
in a domain of ${\bf R}^4$ where $x^1 \neq  0, \;y^1 \neq 0$. Note
that   this example of the web $G_2 (4, 2, 2)$ was also included
into the books  [AG 96], Example {\bf 5.5.2} and [AS 92], Ch.
{\bf 8},  Problem {\bf 5},  and
 the review paper [AG 99a], Section {\bf 6.4},

By (41)--(46), (26), (28), and (31), for the web (50),  we find that
\begin{equation}\label{eq:51}
\renewcommand{\arraystretch}{1.3}
\left\{
\begin{array}{ll}
 a_1 = \displaystyle  \frac{1}{y^1}
- \displaystyle \frac{1}{x^1}, \;\;;\; a_2 = 0; \\
\Gamma_{ij}^1 = \Gamma_{22}^2 = 0, \;\; \Gamma_{12}^2 =-
\displaystyle  \frac{1}{x^1}, \;\; \Gamma_{21}^2 =- \displaystyle
\frac{1}{y^1}, \;\; \Gamma_{11}^2 = \displaystyle  \frac{x^2}{x^1}
- \displaystyle  \frac{y^2}{y^1}, \\ \omega_i^1 = 0, \;\;
\omega_1^2 = \displaystyle
 \Biggl(\frac{x^2}{x^1} - \displaystyle  \frac{y^2}{y^1}\Biggr)
\Biggl(\crazy{\omega}{1}^1 + \crazy{\omega}{2}^1\Biggr), \;\;
\omega_2^2 = -\displaystyle
 \frac{1}{x^1} \crazy{\omega}{1}^1
- \displaystyle  \frac{1}{y^1} \crazy{\omega}{2}^1, \\
      p_{11} = \displaystyle  \frac{1}{(x^1)^2},\;\;
q_{11} = -\displaystyle  \frac{1}{(y^1)^2},\;\; p_{2i} = q_{2i} =
0;
\\
  b_{111}^2 = a_{111}^2 = \Biggl(\displaystyle  \frac{1}{y^1}
- \displaystyle  \frac{1}{x^1}\Biggr) \Biggl(\displaystyle
\frac{x^2}{x^1} - \displaystyle  \frac{y^2}{y^1}\Biggr), \;\;
         b_{121}^2 = - \displaystyle  \frac{1}{(y^1)^2}, \\ b_{112}^2 =
                     \displaystyle  \frac{1}{(x^1)^2},\;\;
          b_{ijk}^1 = b_{211}^2 =
         b_{122}^2 = b_{212}^2 = b_{221}^2 = b_{222}^2 = 0;\\
h_{11} = \displaystyle  \frac{1}{4} \Biggl[\displaystyle
\frac{1}{(y^1)^2} - \displaystyle \frac{1}{(x^1)^2}\Biggr],\;\;
      f_{11} = \displaystyle  \frac{1}{4}
         \Biggl[\displaystyle  \frac{3}{(x^1)^2}
+ \displaystyle  \frac{1}{(y^1)^2}\Biggr], \\
         g_{11} = -\displaystyle  \frac{1}{4} \Biggl[\displaystyle
\frac{1}{(x^1)^2} + \displaystyle  \frac{3}{(y^1)^2}\Biggr], \;\;
         f_{ij} = g_{ij} = h_{ij} = 0, \;\; (i,j) \neq (1,1);
          \\
 - a_{111}^1 = a_{211}^2 = a_{111}^2 =  a_{112}^2
=\displaystyle   \frac{1}{4} \Biggl[\displaystyle
\frac{1}{(x^1)^2} - \displaystyle \frac{1}{(y^1)^2}\Biggr], \\
  a_{ijk}^1 =  a_{122}^2 = a_{212}^2 =
         a_{221}^2 = a_{222}^2  = 0, \;\;
         (i,j,k) \neq (1,1,1) .
        \end{array}
        \right.
\renewcommand{\arraystretch}{1}
\end{equation}

It follows from (51), (25), and (33) that {\em the web $(50)$ is
isoclinic and not transversally geodesic}. Since conditions
(51) imply that
$$
f_{11} + g_{11} + h_{11} =  \frac{1}{4} \Biggl[\displaystyle
\frac{1}{(x^1)^2} - \displaystyle \frac{1}{(y^1)^2}\Biggr],
$$
it follows that {\em the web $(50)$
 belongs to the class {\bf C$_{2}$}}.  By (51), (19), and (22),
 {\it the web $(50)$ belongs to the class {\bf A$_{121}$}}, and {\em it also
belongs to the class {\bf E$_{71}$}}.

{\em The web $(50)$ belongs to the class {\bf F$_2$}} since as was
 proved in  [G 86] (see also [G 88], Example {\bf 8.4.3}),
{\em the three-web $(50)$ can be extended to an exceptional web
$G_2 (4, 2, 2)$ of maximum $2$-rank}. There exists a
two-parameter family of such extensions. The leaves of the fourth
foliation of the web $G_2 (4, 2, 2)$ are defined as level sets of
the  functions
$$
\renewcommand{\arraystretch}{1.3}
\begin{array}{ll}
u_4^1 = \displaystyle \frac{u_3^2 + C_1 u_3^1}{x^2 + y^2 + C_1
u_3^2}, \;\; u_4^2 = - u^1_4 \ln \displaystyle \Biggl|\frac{y^2
+ C_1 y^1 -C_2}{x^2 + C_1 x^1 + C_2}\Biggr| - u_3^1.
\end{array}
\renewcommand{\arraystretch}{1}
$$
The  only abelian equation of the exceptional webs $G_2$ of this
family has the form
\renewcommand{\arraystretch}{1.3}
\begin{eqnarray*}
 - \displaystyle  \frac{1}{x^2 + C_1 x^1 + C_2} dx^1
\wedge dx^2 - \displaystyle    \frac{1}{y^2 + C_1 y^1 - C_2} dy^1
\wedge dy^2  &  & \nonumber \\
 + \displaystyle  \frac{1}{u_3^2 + C_2 u_3^1} du_3^1 \wedge du_3^2
 - \displaystyle  \frac{1}{u^1_4} du_4^1 \wedge du_4^2 = 0. &  &
\end{eqnarray*}
\renewcommand{\arraystretch}{1}
}

\examp{\label{examp:3}
 Consider the web  defined by
\begin{equation}\label{eq:52}
u_3^1 = x^1 + x^2 y^1, \;\; \; u_3^2 =x^2 y^2
\end{equation}
in a domain of ${\bf R}^4$ where $x^2 \neq 0, \; y^2 \neq 0$.

Using (41)--(46), (26), (28), and (31), for the web (52),
we find that
\begin{equation}\label{eq:53}
\renewcommand{\arraystretch}{1.3}
\left\{
\begin{array}{ll}
a_1 = 0, \;\; a_2 =  \displaystyle \frac{1}{x^2 y^2}, \;\;
\Gamma_{21}^1 = \Gamma^2_{22} = -  \displaystyle \frac{1}{x^2
y^2}, \;\; \Gamma_{1k}^i = \Gamma_{21}^2 = \Gamma_{22}^i = 0, \\
\omega_1^1 = -  \displaystyle \frac{1}{x^2 y^2}
\crazy{\omega}{1}^2, \;\; \omega_2^1 = -  \displaystyle
\frac{1}{x^2 y^2} \crazy{\omega}{2}^1, \;\; \omega_1^2 = 0, \;\;
\omega_2^2 = -  \displaystyle \frac{1}{x^2 y^2}
\crazy{\omega}{2}^2, \\ p_{ij} = 0, \;\; (i, j) \neq (2, 2), \;\;
p_{22} = -\displaystyle \frac{1}{(x^2 y^2)^2}, \;\;   q_{ij} = 0,
\\ - b^1_{121} = - b^1_{122} =  b^2_{221} = \displaystyle
\frac{1}{2(x^2 y^2)^2}, \\
 b^1_{111} =  b^2_{222} = b^2_{1kl} =  b^1_{2kl}
= b^2_{212} = b^2_{211} = b^1_{112} = 0, \\
f_{1i} = g_{1i} = h_{1i} = 0, \;\;  g_{22} = h_{22} =
- \displaystyle \frac{1}{2} f_{22} = \displaystyle \frac{1}{3(x^2
y^2)^2}, \\ a^2_{122} = - a^1_{112} = -a^1_{122} =
\displaystyle \frac{1}{6(x^2 y^2)^2}; \; a^i_{jjj} = a^2_{211} =
0.
        \end{array}
        \right.
\renewcommand{\arraystretch}{1}
\end{equation}

It follows from (53), (25), and (33)
 that {\em the web $(52)$ is isoclinic
but not transversally geodesic} since $a^i_{jkl} \neq 0$.
Since
$$
 f_{22} + g_{22} + h_{22} = 0,
$$
{\em the web   $(52)$ belongs to the class {\bf C$_1$}.}
 By (53), (19), and (22),
{\it the web $(52)$ belongs to the class {\bf A$_{121}$}}, and {\em it also
belongs to the class {\bf E$_{2}$}}.

Using the results of  [G 88] (Chapter {\bf 8}), one can prove that
{\em  the three-web $(52)$ cannot be extended to a four-web
$W (4, 2, 2)$ of maximum $2$-rank}. Thus, {\em it belongs to the class {\bf G}}.
}

In Examples {\bf 4--7} below the functions $u_3^1$ and $u_3^2$
will be linear-fractional. In general, a web $W (3, 2, 2)$
of this kind is defined by the equations

\begin{equation}\label{eq:54}
\renewcommand{\arraystretch}{1.3}
\left\{
\begin{array}{ll}
  u_3^1 = \displaystyle \frac{\alpha^1_1 x^1 + \alpha^1_2 x^2 +
  \beta^1_1 y^1 + \beta^1_2 y^2}{\alpha^2_1 x^1 + \alpha^2_2 x^2 +
  \beta^2_1 y^1 + \beta^2_2 y^2}, \\
  {} \\
 u_3^2 = \displaystyle \frac{\alpha^3_1 x^1 + \alpha^3_2 x^2 +
  \beta^3_1 y^1 + \beta^3_2 y^2}{\alpha^4_2 x^1 + \alpha^4_2 x^2 +
  \beta^4_1 y^1 + \beta^4_2 y^2},
\end{array}
\right.
\renewcommand{\arraystretch}{1}
\end{equation}
where the coefficients are arbitrary constants. It is easy to see
that the {\em equations $u_3^1 = \mbox{{\rm const.}}$ and $u_3^2 = \mbox{{\rm const.}} $
define a web consisting of $2$-planes}. The common points
of all 2-planes of this web     are defined by the equations
\begin{equation}\label{eq:55}
\renewcommand{\arraystretch}{1.3}
\left\{
\begin{array}{ll}
\left\{
\begin{array}{ll}
  \alpha^1_1 x^1 + \alpha^1_2 x^2 + \beta^1_1 y^1 + \beta^1_2 y^2 = 0, \\
  \alpha^2_1 x^1 + \alpha^2_2 x^2 + \beta^2_1 y^1 + \beta^2_2 y^2 = 0,
\end{array}
  \right.
  \\
\left\{
\begin{array}{ll}
    \alpha^3_1 x^1 + \alpha^3_2 x^2 + \beta^3_1 y^1 + \beta^3_2 y^2 = 0, \\
  \alpha^4_2 x^1 + \alpha^4_2 x^2 + \beta^4_1 y^1 + \beta^4_2 y^2 = 0.
 \end{array}
 \right.
 \end{array}
 \right.
\renewcommand{\arraystretch}{1}
\end{equation}
Consider the following three matrices of coefficients of system (55):
%\begin{equation}\label{eq:55}
$$
\renewcommand{\arraystretch}{1.3}
D_1 = \pmatrix{\alpha^1_1 & \alpha^1_2 & \beta^1_1 & \beta^1_2, \cr
             \alpha^2_1 & \alpha^2_2 & \beta^2_1 & \beta^2_2   \cr}, \;\;
D_2 = \pmatrix{\alpha^3_1 & \alpha^3_2 & \beta^3_1 & \beta^3_2, \cr
         \alpha^4_2 & \alpha^4_2 & \beta^4_1 & \beta^4_2    \cr},
 \renewcommand{\arraystretch}{1}
$$
and
$$
\renewcommand{\arraystretch}{1.3}
D = \pmatrix{\alpha^1_1 & \alpha^1_2 & \beta^1_1 & \beta^1_2, \cr
             \alpha^2_1 & \alpha^2_2 & \beta^2_1 & \beta^2_2, \cr
             \alpha^3_1 & \alpha^3_2 & \beta^3_1 & \beta^3_2, \cr
        \alpha^4_2 & \alpha^4_2 & \beta^4_1 & \beta^4_2       \cr}.
 \renewcommand{\arraystretch}{1}
$$

%\end{equation}

If  $\mbox{{\rm rank}}\; D_1 = 1$ or  $\mbox{{\rm rank}}\; D_2 =
1$,
then one of the functions $u_3^1$ or $u_3^2$ becomes a constant,
and inequalities in (41) are not satisfied. Thus
equations (54) do not define a web $W (3, 2, 2)$ in these
two cases. Geometrically, in these cases 2-planes defined by (54)
belong to an ${\bf R}^3$.

So, we will assume that $\mbox{{\rm rank}}\; D_1 = \mbox{{\rm rank}}\; D_2 =
2$. The following three cases can occur:
\begin{description}
\item[(i)] $\mbox{{\rm rank}}\; D = 4$. The axes of two pencils of hyperplanes
 of a web defined by (54) have only one common point. Thus
 {\em  the leaves of the foliation $\lambda_3$ are $2$-planes
 passing through a common point but not having a common straight line}.
   This is the general case.

\item[(ii)] $\mbox{{\rm rank}}\; D = 3$. The axes of two pencils of hyperplanes
(and as a result, all 2-planes of the foliation $\lambda_3$) of a
web defined by (54) have a common straight line. Thus in this case
{\em the leaves of  the foliation $\lambda_3$
are $2$-planes of a pencil, and the axis of this pencil is defined by equations
$(55)$}.  This is a special case.

\item[(iii)] $\mbox{{\rm rank}}\; D = 2$.  In this case again
inequalities in (41) are not satisfied.
Geometrically, in this case two pencils of hyperplanes
defined by (54) coincide, their axes coincide, and as a result, all 2-planes
of the foliation $\lambda_3$  defined by (54) coincide with this common
two-dimensional axis. Thus equations (54) do not define a web $W (3, 2, 2)$
in this case.
\end{description}

%If $\det D \neq 0$, then all four hyperplanes defined by each of equations
%(55) are independent. As a result, {\em they and $2$-planes of intersection of
%the first and second pair of hyperplanes have  a unique common
%point, the origin}. This is the general case.

%If $\mbox{{\rm rank}}\; D = 3$, then only three of four hyperplanes defined by each of equa\-tions
%(55) are independent. As a result, {\em they $($and $2$-planes of intersection of
%the first and second pair  of hyperplanes$)$ have  a common straight line
%passing through the origin}.

%Note that if $\mbox{{\rm rank}}\; D = 2$ or $\mbox{{\rm rank}}\; D = 1$, then
%equations  (54) do not satisfy inequalities in (41). Thus
%equations (54) do not define a web $W (3, 2, 2)$ in these
%two cases.
%only two of four hyperplanes defined by each of equations
%(55) are independent. Their intersection is a 2-plane $P$, and {\em all $2$-planes
%of the web in question coincide with  $P$}.

%If $\mbox{{\rm rank}}\; D = 1$, then only one of four hyperplanes defined by each of equations
%(55) is independent. As a result,  {\em all $2$-planes
%of the web in question belongs to that hyperplane}, and the web is
%located in ${\bf R}^3$ instead of  ${\bf R}^4.$ Another way to
%show that equations (54) do not define a web is to check out that
%$\det{(\bar{f}^i_j)} =  \det{(\tilde{f}^i_j)} = 0$,
%and this contradict to the web definition.

The webs $W (3, 2, 2)$ in Examples {\bf 4, 5,} and {\bf 6}
belong to the special case (ii):  2-planes of its foliation
$\lambda_3$ form a pencil.
The web of Example {\bf 7} is of the general type (i):
2-planes of its foliation $\lambda_3$)
have only one common point.

\examp{\label{examp:4}
 Consider the web (see [Bo 35], p. 449) defined by the functions
\begin{equation}\label{eq:56}
u_3^1 =\displaystyle  \frac{y^2}{x^1 + y^1}, \;\;\;
         u_3^2 =\displaystyle  \frac{x^2}{x^1  + y^1}
\end{equation}
in a domain of ${\bf R}^4$ where $x^2 \neq 0, \; y^2 \neq 0, \;
x^1 + y^1 \neq 0$.

In this case the leaves of  the foliations $\lambda_1$ and
$\lambda_2$  are  2-planes parallel to
the coordinate 2-planes $O y^1 y^2$ and $O x^1 x^2$, respectively.
 Since ${{\rm rank}}\; D_1 = {{\rm rank}}\; D_2 = 2, \;{{\rm rank}}\; D = 3$,
 the leaves of the third foliation $\lambda_3$ are  2-planes  of a pencil.
  It is easy to see that  the axis of this pencil is
 defined by the equations $x^1 + y^1 = 0, \; x^2 = 0, \; y^2 = 0$.

Using (41)--(46), (26), (28), and (31), for the web (56),
 we find that
\begin{equation}\label{eq:57}
\renewcommand{\arraystretch}{1.3}
\left\{
\begin{array}{ll}
\Gamma^1_{2k} = \Gamma^2_{k1} = 0, \;\; \Gamma^1_{11} =
\Gamma^2_{12} = -\sigma, \;\;
\Gamma^1_{12} = \Gamma^2_{22} = -\tau, \\
\omega_1^1 = - [\sigma(\crazy{\omega}{1}^1  + \crazy{\omega}{2}^2)
+ \tau \crazy{\omega}{2}^1], \;\; \omega_1^2 = - \sigma \crazy{\omega}{2}^2, \\
\omega_2^2 = - [\tau (\crazy{\omega}{1}^2  + \crazy{\omega}{2}^2)
+ \sigma \crazy{\omega}{1}^1], \;\; \omega_2^1 = - \tau \crazy{\omega}{1}^1, \\
a_1 = -\sigma, \;\; a_2 =\tau,  \;\; b^i_{jkl} = a^i_{jkl} =
0, \;\; p_{ij} = q_{ij} = f_{ij} = g_{ij} =  h_{ij} = 0,
\end{array}
\right.
\renewcommand{\arraystretch}{1}
\end{equation}
where $ \sigma = \displaystyle  \frac{x^1 + y^1}{y^2}$
and $\tau =  \displaystyle  \frac{x^1 + y^1}{x^2}$.

It follows from (57), (32), and (39)
 that {\em the web $(56)$ is an algebraizable
 group nonparallelizable web}, i.e., {\em it belongs to the class
{\bf  D$_{231}$}}. By (57) and (8), {\em it belongs also
to the classes {\bf A$_{1}$} and {\bf E$_{1}$}}.

%0ne can also prove (see [Bo 35], p. 449)
%that {\em this web is equivalent to a web formed by three pencils
%of $2$-planes in ${\bf R}^4$ whose axes do not belong to an ${\bf
%R}^3$.}

It is well-known that in a representation of
 a Grassmannizable (transversally geodesic and
isoclinic) web $W (3, 2, 2)$, such a web
 is generated by  three surfaces $X_\xi, \; \xi = 1, 2, 3,$ of general
position in the projective space $P^3$, and for an  algebraizable
(Grassmannizable and hexagonal), these three surfaces
 belong to a cubic surface $V^2_3 \subset P^3$
(see [A 73] or [AS 92], Section {\bf 3.3}). The tensors $f_{ij},
g_{ij}$, and $ h_{ij}$ are the second fundamental tensors of the
surfaces $X_\xi$. The conditions $f_{ij} = g_{ij} =
h_{ij} = 0$ imply that  all three surfaces $X_\xi$ are planes.

Using the results of  [G 88],  Example {\bf
8.4.1}, one can prove that {\em  the three-web $(56)$ can be
extended to a four-web $W (4, 2, 2)$ of maximum $2$-rank, and the
latter is generated in  $P^3$ by four planes of general position}.
  Thus, {\em the web $(56)$ belongs to the class {\bf  F$_{1}$}}.
}

\examp{\label{examp:5}
 Consider the web defined by the functions
\begin{equation}\label{eq:58}
u_3^1 =\displaystyle  \frac{x^1 + y^1}{x^1 + x^2}, \;\;\; u_3^2
=\displaystyle  \frac{x^2 + y^2}{y^1 + y^2}
\end{equation}
in a domain of ${\bf R}^4$, where $x^1 + x^2 \neq 0,\; y^1 + y^2
\neq 0, \; x^2 \neq y^1$.
%the leaves of
% the first (or the second) foliation of this web are pencils of
% parallel 2-planes,  while

 In this case the leaves the foliations $\lambda_1$ and $\lambda_2$   are
 2-planes parallel to
the coordinate 2-planes $O y^1 y^2$ and $O x^1 x^2$, respectively.
 Since ${{\rm rank}}\; D_1 = {{\rm rank}}\; D_2 = 2, \;{{\rm rank}}\; D  = 3$,
 the leaves
 of the third foliation $\lambda_3$ are  2-planes  of a pencil.
 It is easy to see that  the axis of this pencil is
 defined by the equations  $x^1 = y^2 = - x^2 = -y^1$.

Using (41)--(46), (26), (28), and (31), for the web (58),
 we find that
\begin{equation}\label{eq:59}
\renewcommand{\arraystretch}{1.3}
\left\{
\begin{array}{ll}
\Gamma^1_{11} = - \Gamma^2_{21} = \sigma, \;\; \Gamma^1_{21} = -
\Gamma^2_{22} = \tau, \;\; \Gamma^1_{12} = \Gamma^1_{22} =
\Gamma^2_{11} = \Gamma^2_{21} =  0, \\ \omega_1^1 = \sigma
(\crazy{\omega}{1}^1 + \crazy{\omega}{2}^1) + \tau
\crazy{\omega}{1}^2, \;\; \omega_1^2 = - \sigma
\crazy{\omega}{1}^2, \\ \omega_2^2 =  - \tau (\crazy{\omega}{1}^2
+ \crazy{\omega}{2}^2) - \sigma \crazy{\omega}{2}^2, \;\;
\omega_2^1 =  \tau   \crazy{\omega}{2}^1, \\ a_1 = \sigma,  \;\;
a_2 = \tau,   \;\; p_{ij} =  q_{ij} = 0, \\ b^i_{jkl} = a^i_{jkl}
= 0, \;\; f_{ij} = g_{ij} = h_{ij} = 0.
        \end{array}
        \right.
\renewcommand{\arraystretch}{1}
\end{equation}
where $\sigma = \displaystyle  \frac{x^1 + x^2}{x^2 - y^1}$ and
$\tau =  \displaystyle \frac{y^1 + y^2}{x^2 - y^1}$.

It follows from (59), (32), and (39)
 that {\em the
 web  $(58)$ belongs to the same classes {\bf A$_{1}$}, {\bf
 D$_{231}$},  {\bf   E$_{1}$}, and {\bf  F$_{1}$} as the web $(56)$}.

 %Since
 %geometrically the webs (72) and (74) are the same, they are equivalent, and
 %for the web (74) we have equations similar to equations (73), and this
% web belongs to the same classes {\bf A${}_{1}$},
% {\bf  D${}_{231}$}, and {\bf  F${}_{1}$}. Note that the equivalence of
% webs (72) and (74) can be established by the following transformations
%of variables $x^i$ and the variables $y^i$:
%\begin{equation}\label{eq:75}
%x^2 \rightarrow (x^1 + x^2), \;\; y^2 \rightarrow (y^1 + y^2),
%\end{equation}
%transferring the axis of the third pencil of the web (72) to
%the axis of the third pencil of the web (72) and not changing
%the first two foliations of webs.
}

\examp{\label{examp:6}
 Consider the web defined by the functions
\begin{equation}\label{eq:60}
u_3^1 =\displaystyle  \frac{x^1 + y^1}{x^1 + x^2}, \;\;\; u_3^2
=\displaystyle  \frac{x^2 + y^2}{x^1 + x^2}
\end{equation}
in a domain of ${\bf R}^4$ where $x^1 + x^2 \neq 0,\; y^1 + y^2
\neq 0$. Chern in [C 36], pp. 357--358, gave an example of a web
$W (3, 2, r)$ formed in ${\bf R}^{2r}$ by three pencils of $r$-planes with fixed
$(r-1)$-dimensional axes that do not belong to any
hyperplane ${\bf R}^{2r - 1}$.
The web (60) is the Chern's web corresponding to the
value $ r = 2$.
%Note that it is different from the webs (72) and (74).
%the leaves of
% the first (or the second) foliation of this web are pencils of
% parallel 2-planes, while

 Since ${{\rm rank}}\; D_1 = {{\rm rank}}\; D_2 = 2, \;{{\rm rank}}\; D = 3$,
 the leaves
 of the third foliation $\lambda_3$  are  2-planes  of a pencil.
  It is easy to see that  the axis of this pencil is
 defined by the equations $x^1 = y^2 = - x^2 = -y^1$.

 Using (41)--(46), (26), (28), and (31), for the web (60),
 we find that
\begin{equation}\label{eq:61}
\renewcommand{\arraystretch}{1.3}
\left\{
\begin{array}{ll}
\Gamma^1_{11} = \Gamma^1_{21} = \Gamma^2_{12} = \Gamma^2_{22} = -
\sigma, \;\; \Gamma^1_{12}= \Gamma^1_{22} = \Gamma^2_{11} =
\Gamma^2_{21} =0, \\
a_1 = a_2 =  - \sigma, \;\;        p_{ij} =  q_{ij} = 0,  \\
\omega_1^1 = - \sigma (\crazy{\omega}{1}^1 + \crazy{\omega}{2}^1
              + \crazy{\omega}{1}^2), \;\; \omega_1^2 = - \sigma \crazy{\omega}{2}^2. \\
\omega_2^2 = - \sigma (\crazy{\omega}{1}^1 + \crazy{\omega}{2}^2
               + \crazy{\omega}{1}^2), \;\;
\omega_2^1 =  -\sigma \crazy{\omega}{2}^1. \\ \displaystyle
\frac{1}{2} b^1_{111} = \displaystyle \frac{1}{2} b^2_{222} =
\displaystyle \frac{1}{2} b^1_{121} =
  \displaystyle \frac{1}{2} b^2_{212} =  b^2_{122}
=  b^1_{112} = b^1_{122} = b^1_{211} \\= b^1_{221} = b^2_{211} =
b^2_{112}  = b^2_{221}
 = \sigma^2, \;\; b^1_{222} = b^2_{111} = b^1_{212} = b^2_{121}
= 0, \\ a^1_{111} = a^2_{222}
 = 2 a^1_{112} =  2 a^2_{122} =  \displaystyle \frac{4}{3}\sigma^2, \;\;
 a^1_{122} =  a^2_{112} = a^1_{222} = a^2_{111} = 0, \\
 f_{ij} = g_{ij} = h_{ij} = \displaystyle \frac{2}{3}\sigma^2,
        \end{array}
        \right.
\renewcommand{\arraystretch}{1}
\end{equation}
where $\sigma = \displaystyle \frac{x^1 + x^2}{y^1 + y^2}$.

It follows from (61), (25), and (33) that {\em the web $(60)$ is isoclinic
but not transversally geodesic} since, for example, $a^1_{111} =
 \displaystyle \frac{4}{3} \sigma^2 \neq 0$. Thus, {\em this web
 belongs to the class {\bf C}}.
  By (61), (21),  and (24), {\em it belongs also
to the classes {\bf E$_{1}$} and {\bf A$_{1121}$}}.

 Using the results of Ch. 8 of [G 88], one can prove that
{\em the web $(60)$ cannot be extended to a web $W (4, 2, 2)$ of
maximum $2$-rank}. Thus, {\em it belongs to the class {\bf G}}.
}

In the next  example the third foliation of a web $W (3, 2, 2)$
is also formed by 2-planes. However, unlike the webs
in Examples {\bf 4, 5}, and {\bf 6}, these 2-planes do not belong to
a pencil---they have only one common point.

\examp{\label{examp:7}
 Consider the
web (see [Bo 35], p. 449; [C 36], p. 357) defined by
\begin{equation}\label{eq:62}
u_3^1 =\displaystyle  \frac{x^1 + y^1}{x^2 + y^2}, \;\;\; u_3^2
=\displaystyle  \frac{x^1 - y^1}{x^2 - y^2}
\end{equation}
in a domain of ${\bf R}^4$ where $x^2 \neq \pm y^2, \; x^2 y^1 -
x^1 y^2 \neq 0$.

Since in this case ${{\rm rank}}\; D_1 = {{\rm rank}}\; D_2 = 2, \;
{{\rm rank}}\; D = 4$, the leaves
 of the third foliation  $\lambda_3$ are  2-planes  not belonging to a pencil.

Using (41)--(46), (26), (28), and (31), for the web (62),
we find that
\begin{equation}\label{eq:63}
\renewcommand{\arraystretch}{1.3}
\left\{
\begin{array}{ll}
\Gamma^1_{11} = -\Gamma^2_{22} = - \rho, \;\; \Gamma^1_{12} =
-\Gamma^1_{21} = - \sigma, \\ \Gamma^2_{12} =  - \Gamma^2_{21} = -
\tau, \;\; \Gamma^1_{22} = \Gamma^2_{11} = 0,\\ \omega_1^1 = -
\rho  (\crazy{\omega}{1}^1  + \crazy{\omega}{2}^1) + \sigma
(\crazy{\omega}{1}^2 -   \crazy{\omega}{2}^2), \;\; \omega_1^2 =
\tau (\crazy{\omega}{1}^2 -  \crazy{\omega}{2}^2),
\\
\omega_2^2 = - \tau  (\crazy{\omega}{1}^1 - \crazy{\omega}{2}^1) +
\rho (\crazy{\omega}{1}^2 +  \crazy{\omega}{2}^2), \;\; \omega_2^1
= - \displaystyle  \frac{x^1 + y^1}{y^2} \crazy{\omega}{1}^1, \\
a_1 =  - 2\tau, \;\;  a_2 = 2\sigma, \;\; p_{11} = - q_{11}
 = -2 \tau^2,
\\ p_{22} = - q_{22}
=- 2 \sigma^2, \;\; p_{12} = p_{21} = - q_{12} = - q_{21} = -
\displaystyle \frac{1}{2} \rho^2, \\ b^1_{221} = - b^1_{212}=
p_{22}, \;\; b^1_{121} = - b^1_{112}= b^2_{212} = - b^2_{221}=
p_{12}, \\ b^2_{112} = - b^2_{121}= p_{11}, \;\; b^i_{jjj} =
b^i_{211} = b^i_{122} =  0, \\ f_{ij} = - g_{ij} = p_{ij}, \;\;
h_{ij} = 0,\;\; b^i_{(jkl)} = a^i_{jkl} = 0.
        \end{array}
        \right.
\renewcommand{\arraystretch}{1}
\end{equation}
where
$$
\sigma \displaystyle = \frac{(x^2 - y^2)^2}{2(x^1 y^2 -
x^2 y^1)}, \;\; \tau = \displaystyle \frac{(x^2 + y^2)^2}{2(x^1
y^2 - x^2 y^1)}, \;\; \rho = \displaystyle \frac{(x^2)^2 -
(y^2)^2}{x^1 y^2 - x^2 y^1},
$$
and $ \rho^2 = 4 \sigma \tau$.

By (63), equation (8) does not hold.
Thus {\em the web $(62)$ belongs to the class {\bf A${}_{2}$}}.
It follows from (63), (33), and
(38) that {\em the web $(62)$ is algebraizable Bol web $B_m$ but
not a group web} since $b^i_{jkl} \neq 0$. By (63) and (38),
{\em this web belongs  to the class {\bf D${}_{21}$}.
It also belongs to the class {\bf E${}_{8}$}}.

One can prove (see [Bo 35], p. 448) that this web can be realized as follows. Take four
generators $l_1, l_2, l_3$ and $l_4$ of a quadric $Q \subset {\bf
R}^4$ in such a way that $(l_1, l_2, l_3, l_4) = -1$. The quadric
$Q$ always belong to some ${\bf R}^3$. Take a point $P \notin {\bf
R}^3$.
 The straight
lines intersecting $l_3$ and $l_4$ form a linear congruence $L$ in
${\bf R}^3$. The foliations $\lambda_\xi$ of the web (62) are:
$\lambda_1$ consists of 2-planes of the pencil with the axis
$l_1$; $\lambda_2$ consists of 2-planes of the pencil with the
axis $l_2$; and $\lambda_3$ consists of 2-planes passing through
the point $P$ and rays of $L$.

A realization of the web (62) in $P^3$ is generated  by a cubic surface
 $V^2_3$ which decomposes into a quadric $V^2_2$ and a plane $V^2_1$.
Using the results of  [G 88] (Chapter {\bf 8},
Example {\bf 8.4.1}), one can prove that {\em  the three-web
$(62)$ can be extended to a four-web $W (4, 2, 2)$ of maximum
$2$-rank, and the latter is generated in $P^3$ by a quartic which decomposes
into a quadric $V^2_2$ and two planes  $V^2_1$ and
$\widetilde{V}^2_1$}. This four-web $W (4, 2, 2)$
is not exceptional since it is algebraizable. Thus {\em the web
$(62)$ belongs to the class  {\bf F${}_{1}$}}.
}

\examp{\label{examp:8}
 Consider the web  defined by the functions
\begin{equation}\label{eq:64}
u_3^1 = A x^1 y^1 + B x^1 y^2 +C x^2 y^1 + E x^2 y^2, \;\;\; u_3^2
= a x^1 y^1 + b x^1 y^2 +b x^2 y^1 + e x^2 y^2
\end{equation}
(where $A, B, C, E, a, b, c$, and $e$ are constants) in a domain
of ${\bf R}^4$ where
$$
\renewcommand{\arraystretch}{1.3}
\left\{
\begin{array}{ll}
\Delta_1 = (y^1)^2 (A c - Ca) + (y^2)^2 (B e - E b)
+ y^1 y^2 (A e + B c - C b - E a) \neq 0, \\
\Delta_2 = (x^1)^2 (A b - B a) + (x^2)^2 (C e - E c)
+ x^1 x^2 (A e - B c + C - E a) \neq 0.
        \end{array}
        \right.
\renewcommand{\arraystretch}{1}
$$
It is naturally to call webs (64) {\em homogeneous quadratic webs}.

Using (41)--(46), (26), (28), and (31), for the web (64),
we find that for arbitrary  $A, B, C, E, a, b, c$, and $e$  we have
\begin{equation}\label{eq:65}
a_1 = a_2 = 0, \;\; p_{ij} = q_{ij} =  0, \;\;
 f_{ij} = g_{ij} = h_{ij} = 0, \;\; b^1_{111} = a^1_{111} \neq 0.
\end{equation}

To prove the last inequality, we calculate  the  term of
$b^1_{111}$ with $(x^1)^2 (y)^1)^2$. Rather difficult calculations give
that the coefficient in this term is
$$
 a \Bigl[ (A b - B a) (- 2 B c^2 + 2 E  a c - C a e + A c e)
  - (C b - E a)  (A c - C a)  b\Bigl].
$$
This shows that for  arbitrary  $A, B, C, E, a, b, c$, and $e$
we have $ b^1_{111} = a^1_{111} \neq 0$.
 However, for some values of these
 constants, we have $b^i_{jkl} = a^i_{jkl} = 0$ for all components
 (see Example {\bf 9}, where
$B = C = a = e =0, \; A = E = b = c = 1$).

It follows from (65), (32), and (33)
 that {\em the web $(64)$ is always isoclinicly geodesic and
in general not transversally geodesic}.
Thus, {\em this web belongs to the classes {\bf B} and {\bf C}}.
Conditions (65) show that {\em it also belongs to the class ${\bf E_1}$}.

It is an interesting problem to {\em determine all homogeneous quadratic
group $($and therefore parallelizable$)$  webs}
(i.e., webs (64) with the vanishing curvature tensor).
The web in our next example is of this kind.
}

\examp{\label{examp:9}
 Consider the web  defined by
\begin{equation}\label{eq:66}
u_3^1 = x^1 y^1 + x^2 y^2, \;\; \; u_3^2 = x^1 y^2 + x^2 y^1
\end{equation}
in a domain of ${\bf R}^4$ where $x^2 \neq \pm x^1, \; y^2 \neq
\pm y^1$. Note that the web (66) is a particular case of the web (64) corresponding
to the following values of constants: $B = C = a = e =0, \;
A = E = b = c = 1$.

Using (41)--(46), (26), (28), and (31), for the web (66),
we find that
\begin{equation}\label{eq:67}
\renewcommand{\arraystretch}{1.3}
\left\{
\begin{array}{ll}
\Gamma^1_{11} = \Gamma^1_{22} = \Gamma^2_{12} = \Gamma^2_{21} = -
\displaystyle \frac{x^1 y^1 + x^2 y^2}{\bigl((x^1)^2 -
(x^2)^2\bigr) \bigl((y^1)^2 -  (y^2)^2\bigr)},
\\
\Gamma^1_{12} = \Gamma^1_{21} = \Gamma^2_{11} = \Gamma^2_{22} =
\displaystyle \frac{x^1 y^2 + x^2 y^1}{\bigl((x^1)^2 -
x^2)^2\bigr)
       \bigl((y^1)^2 -  (y^2)^2\bigr)},\\
\omega_j^i = \Gamma^i_{jk} (\crazy{\omega}{1}^k +
\crazy{\omega}{2}^k),
%        + \Gamma_{12}^1 (\crazy{\omega}{1}^2 + \crazy{\omega}{2}^2), \;\;
%\omega_1^2 = \Gamma^2_{12} (\crazy{\omega}{1}^2 +
%\crazy{\omega}{2}^2)
%        + \Gamma_{11}^2 (\crazy{\omega}{1}^1 + \crazy{\omega}{2}^1), \\
%\omega_2^2 =  \Gamma^2_{12} (\crazy{\omega}{1}^1 +
%\crazy{\omega}{2}^1)
%       + \Gamma_{12}^2 (\crazy{\omega}{1}^2 + \crazy{\omega}{2}^2), \;\;
%\omega_2^1 =  \Gamma^1_{22} (\crazy{\omega}{1}^2 +
%\crazy{\omega}{2}^2)
%      + \Gamma_{12}^1 (\crazy{\omega}{1}^2 + \crazy{\omega}{2}^1), \;\;
\;\; a_1 = a_2 = 0, \;\; p_{ij} =  q_{ij} = 0,  \\ b^i_{jkl} =
a^i_{jkl} = 0,  \;\; f_{ij} = g_{ij} = h_{ij} = 0.
        \end{array}
        \right.
\renewcommand{\arraystretch}{1}
\end{equation}

It follows from (67), (32), and (39)
 that {\em the web $(66)$ is parallelizable}.
 Thus, {\it this web belongs
to the class {\bf D${}_{232}$} $($and therefore to the class
{\bf B}$)$}.
 By (67) and (8), {\em it belongs also
to the classes {\bf E${}_1$} and {\bf A${}_1$}}.

Since this web is  algebraizable,
its representation  is generated by three planes
of a pencil of planes
in $P^3$. Using the results of  [G 88] (Chapter {\bf 8},
Example {\bf 8.4.1}), one can prove that {\em  the three-web
$(66)$ can be extended to a four-web $W (4, 2, 2)$ of maximum
$2$-rank, and the latter is generated by four planes belonging to
a pencil of planes in $P^3$.}  This four-web $W (4, 2, 2)$
is not exceptional since it is algebraizable. Thus {\em the web
$(66)$ belongs to the class  {\bf F${}_{1}$}}.
}

\examp{\label{examp:10}
 Consider the web  defined by
\begin{equation}\label{eq:68}
u_3^1 = x^1 e^{-y^2} + x^2 e^{y_1}, \;\;\; u_3^2 =  x^2 + y^2
\end{equation}
in a domain of ${\bf R}^4$ where $x^2 \neq 0$.
Using (41)--(46), (26), (28), and (31), for the web (68),
we find that
\begin{equation}\label{eq:69}
\renewcommand{\arraystretch}{1.3}
\left\{
\begin{array}{ll}
\Gamma^2_{ij} = \Gamma^1_{11} = 0, \;\; \Gamma^1_{12} = 1, \;\;
\Gamma^1_{21} = - \displaystyle \frac{1}{x^2}, \;\; \Gamma^1_{22} =
- \Bigl(e^{y^1} +  \displaystyle \frac{x^1}{x^2} e^{-y^2}\Bigr),
\\ \omega_i^2 = 0, \;\; \omega_1^1 =  \displaystyle \frac{1}{x^2}
\crazy{\omega}{1}^2, \;\; \omega_2^1 =  \displaystyle
\frac{1}{x^2} \crazy{\omega}{1}^1 - \Bigl(e^{y^1} +  \displaystyle
\frac{x^1}{x^2} e^{-y^2}\Bigr), \;\;
a_1 = 0, \\ a_2 = - \Biggl(1
+ \displaystyle \frac{1}{x^2}\Biggr),  \;\; p_{ij} = 0, \;\; (i, j)
\neq (2, 2),\;\; p_{22} = \displaystyle \frac{1}{(x^2)^2}, \;\;
q_{ij} = 0, \\ b^1_{222} = a_{222}^1 = \displaystyle \frac{1}{x^2}
\Biggl(e^{y^1} +  \frac{x^1}{x^2} e^{-y^2}\Biggr), \;\; b_{221}^1
=  \frac{1}{(x^2)^2}, \\ b^1_{jkl} = 0, \; (j, k, l) \neq (2,2,2),
(2,2,1), \;\; b^2_{jkl} = 0, \\
a^1_{122} =
=  \displaystyle \frac{5}{24 (x^2)^2}, \;\; a^2_{222} = -
\frac{3}{8 (x^2)^2}, \;\; a_{111}^i = a_{112}^i = 0, \\
f_{1i} = g_{1i} = h_{1i} = 0, \; f_{22} = \displaystyle \frac{19}{24 (x^2)^2},
 \; \; g_{22} = h_{22} = -\displaystyle \frac{5}{24 (x^2)^2}.
        \end{array}
        \right.
\renewcommand{\arraystretch}{1}
\end{equation}

It follows from (69), (25),  and (33)
 that {\em the web $(68)$ is isoclinic
but  not transversally geodesic}.
Thus, {\em this web belongs to the class {\bf C}}.
 By (69), (20), and (23), {\em this web belongs to the
 class {\bf A$_{131}$};
 it also belongs to the class {\bf E$_2$}}.

Using the results of Ch. 8 of [G 88], one can prove that {\em
 the web $(68)$ cannot be extended to a web $W (4, 2, 2)$
 of maximum $2$-rank}.
Thus, {\em it belongs to the class {\bf G}}.
}

\examp{\label{examp:11}
 Consider the web defined by
\begin{equation}\label{eq:70}
u_3^1 = e^{-2 y^2} \bigl(x^1 + e^{2y^2} + \displaystyle
\frac{e}{2} x^2 y^2\bigr), \;\;\; u_3^2 =x^2 + y^2
\end{equation}
in any domain of ${\bf R}^4$.

Using (41)--(46), (26), (28), and (31), for the web (70),
we find that for this web
\begin{equation}\label{eq:71}
\renewcommand{\arraystretch}{1.3}
\left\{
\begin{array}{ll}
\Gamma^2_{ij} = \Gamma^1_{11} = 0, \;\; \Gamma^1_{12} = 2, \;\;
\Gamma^1_{21} = \Bigl(y^2 - \displaystyle \frac{1}{2}\Bigr)
e^{1-2y^2},\\
\Gamma^1_{22} = \Bigl[-y^2 + e^{-2y^2} \Bigl(y^2
- \displaystyle \frac{1}{2}\Bigr)\Bigl(2x^1 + e x^2 y^2
- \displaystyle \frac{1}{2} e y^2\Bigr)\Bigr] e^{1-2y^2},\\
\omega_1^1
= \Bigl(y^2 - \displaystyle \frac{1}{2}\Bigr) e^{1-2y^2}
\crazy{\omega}{1}^2 + 2 \crazy{\omega}{2}^2, \;\; \omega_i^2 = 0,\\
\omega_2^1 = 2 \crazy{\omega}{1}^1 +  \Bigl(y^2 - \displaystyle
\frac{1}{2}\Bigr) e^{1-2y^2} \crazy{\omega}{2}^1 \\
+ \Bigl[-y^2 +
e^{-2y^2} \Bigl(y^2 - \displaystyle \frac{1}{2}\Bigr)\Bigl(2x^1 + e x^2 y^2
- \displaystyle \frac{1}{2} e y^2\Bigr)\Bigr]
e^{1-2y^2}  (\crazy{\omega}{1}^2 + \crazy{\omega}{2}^2), \\
 a_1 = 0, \;\;
a_2 = \Bigl(y^2 - \displaystyle \frac{1}{2}\Bigr) e^{1-2y^2} - 2,\\
p_{1i} =  q_{1i} = p_{21} = p_{22} = q_{21} = 0, \;\; q_{22} =
2 (1 - y^2)e^{1-2y^2},  \\ b^2_{jkl} = b^1_{111} = b^1_{112} =
b^1_{121} = b^1_{211} =0, \\ b^1_{222} = (2y^2 - 1)\Bigl[- \Bigl(1
+ \displaystyle \frac{1}{2} ey^2\Bigr) e^{1-2y^2} + \Bigl(y^2
-\displaystyle \frac{1}{2}\Bigr) e^{2 - 3y^2}\\
+ \Bigl(4 x^1 + 2 e x^2 y^2 - \displaystyle \frac{1}{2}
ex^2\Bigr) e^{1-4y^2}\Bigr],\;\; b^1_{221} = \displaystyle
\frac{1}{2} e^{1-2y^2},\\ b^1_{122} = \Bigl(y^2 - \displaystyle
\frac{1}{2}\Bigr) e + (y^2 - 1) e^{1-2y^2}, \;\; b^1_{212}=
\Bigl(y^2 - \displaystyle \frac{1}{2}\Bigr) (e -  e^{1-2y^2}),\\
f_{1i} = g_{1i} = h_{1i} = 0, \\
 f_{22} = h_{22} = \displaystyle \frac{3e}{16} \Bigl(2 y^2 - 1\Bigr) +
\Bigl(\displaystyle \frac{2}{3} y^2  - \displaystyle
\frac{19}{24}\Bigr) e^{1-2y^2}, \\ g_{22} = \displaystyle
\frac{3e}{16} \Bigl(2 y^2 - 1\Bigr) + \Bigl(\displaystyle
\frac{29}{24}  - \displaystyle \frac{4}{3}y^2\Bigr) e^{1-2y^2}, \\
%f_{22} + g_{22} + h_{22} = \displaystyle \frac{3e}{8} \Bigl(2 y^2
%- 1\Bigr) + \Bigl(\displaystyle \frac{5}{12} - \displaystyle
%\frac{4}{3} y^2\Bigr) e^{1-2y^2}, \\
a^2_{jkl} = a^1_{111} =
a^1_{112} = 0, \;\;
 a^1_{222} = b^1_{222}, \\
 a^1_{122} = \displaystyle \frac{1}{8}
\Bigl(2 y^2 - 1\Bigr) \Bigl(2e - 1\Bigr) + \displaystyle
\frac{1}{6} \Bigl(\displaystyle \frac{4}{9} y^2 - \displaystyle
\frac{41}{36}\Bigr) e^{1-2y^2}.
        \end{array}
        \right.
\renewcommand{\arraystretch}{1}
\end{equation}

It follows from (71), (25),  and (33) that {\em the web $(70)$ is
isoclinic but not transversally geodesic}.
Thus, {\em this web belongs to the class {\bf C}}.
 By (71), (20), and (23), {\em this web belongs to the class {\bf A$_{131}$};
 it also belongs to the class ${\bf E}_3$}.
}

\examp{\label{examp:12}
 Consider the web  defined by the functions
\begin{equation}\label{eq:72}
u_3^1 =   \displaystyle \frac{(x^1)^2 + (y^1)^2}{x^1 +  y^1},
\;\;\; u_3^2 =  \frac{(x^2)^2 + (y^2)^2}{x^2 +  y^2},
\end{equation}
 in a domain of ${\bf R}^4$ where $x^1 \neq - y^1, \; x^2 \neq y^2, \;
 (x^i)^2 + 2 x^i y^i - (y^i)^2   \neq 0, \;\;
  (y^i)^2 + 2 x^i y^i - (x^1)^2   \neq 0$.

Using (41)--(46), (26), (28), and (31), for the web (72),
we find that for this web
\begin{equation}\label{eq:73}
\renewcommand{\arraystretch}{1.3}
\left\{
\begin{array}{ll}
\Gamma^1_{11} = \displaystyle \frac{4 x^1 y^1 (x^1 +
y^1)}{(A^1)^2}, \;\;\;
 \Gamma^1_{jk} = 0, \;\; (j, k) \neq (1, 1), \\
\Gamma^2_{22} = \displaystyle \frac{4 x^2 y^2 (x^2 +
y^2)}{(B^2)^2}, \;\;\; \Gamma^2_{jk} = 0, \;\; (j, k) \neq (2, 2), \\
\omega_j^i = 0, \;\; j \neq i; \;\; \omega_1^1 =
\displaystyle \frac{4 x^1 y^1 (x^1 + y^1)}{(A^1)^2}
\crazy{\omega}{1}^1, \;\;\; \omega_2^2 = \displaystyle \frac{4 x^2
y^2 (x^2 + y^2)}{(B^2)^2} \crazy{\omega}{2}^2, \\ a_1 = 0, \;\;\;
a_2 = 0,\;\;\; p_{ij} = q_{ij} = 0, \\
b^1_{111}  = -\displaystyle
\frac{4 (x^1 + y^1)^2 [(x^1)^4 - 2 (x^1)^3 y^1 +  6 (x^1)^2
(y^1)^2 + 2 x^1 (y^1)^3 + (y^1)^4]}{(A^1)^4 A^2}, \\ b^2_{222}  =
\displaystyle \frac{4 (x^2 + y^2)^2 [(x^2)^4 + 2 (x^2)^3 y^2 + 6
(x^2)^2 (y^2)^2 - 2 x^2 (y^2)^3 + (y^2)^4]}{(B^1)^4 B^2}, \\
b^1_{jkl} = 0, \;\; (j, k, l) \neq (1, 1, 1), \;\; b^2_{jkl} = 0,
\;\; (j, k, l) \neq (2, 2, 2), \\ a^i_{jjj} = a^i_{iij} = 0, \;\;
i \neq j, \;\; a^1_{111} = \displaystyle \frac{3}{4} b^1_{111},
\;\; a^1_{122} = -\displaystyle \frac{1}{12} b^2_{222}, \\
a^2_{222} = \displaystyle \frac{3}{4} b^2_{222}, \;\; a^2_{112} =
-\displaystyle \frac{1}{12} b^1_{111}, \\
f_{11} = g_{11} = h_{11}
= \displaystyle \frac {1}{4} b^1_{111}, \;\; f_{12} = g_{12} =
h_{12} = 0, \;\; f_{22} = g_{22} = h_{22} = \displaystyle \frac
{1}{4} b^2_{222}.
        \end{array}
        \right.
\renewcommand{\arraystretch}{1}
\end{equation}
where $$ A^i = (x^i)^2 + 2 x^i y^i - (y^i)^2, \;\; B^i = (y^i)^2 +
2 x^i y^i - (x^i)^2. $$

It follows from (73), (32), and (33)
 that {\em the web $(72)$ is  isoclinicly geodesic
and not transversally geodesic}. Thus, {\em this web belongs to
the classes {\bf B} and {\bf C}}. Conditions (73) show that {\em
it also belongs to the class {\bf E$_1$}}.
}

\examp{\label{examp:13}
 Consider the web  defined by the functions
\begin{equation}\label{eq:74}
u_3^1 =  x^1 +  y^1 +  \displaystyle \frac{1}{2} (x^1)^2 y^1,
\;\;\; u_3^2 =  x^2 + y^2 +  \displaystyle \frac{1}{2} (x^2)^2 y^2
\end{equation}
 in a domain of ${\bf R}^4$ where $x^1 y^1 \neq -1, \; x^2 y^2  \neq -1$.

Using (41)--(46), (26), (28), and (31), for the web (74),
we find that
%\newpage
\begin{equation}\label{eq:75}
\renewcommand{\arraystretch}{1.5}
\left\{
\begin{array}{ll}
\Gamma^1_{i2} =   \Gamma^2_{1i} = \Gamma^1_{21} =
\Gamma^2_{21} = 0, \;\; \Gamma^1_{11} =
-\displaystyle \frac{x^1}{A^1 (B^1)^2}, \;\;  \;\; \Gamma^2_{22} = -\displaystyle
\frac{x^2}{A^2 (B^2)^2}, \\ \omega_j^i = 0, \;\; j \neq i; \;\;
\omega_1^1 =   -\displaystyle \frac{x^1}{A^1 (B^1)^2}
\crazy{\omega}{1}^1, \;\;\; \omega_2^2 = -\displaystyle
\frac{x^2}{A^2 (B^2)^2} \crazy{\omega}{2}^2, \\ a_1 = 0, \;\;\;
a_2 = 0,\;\;\; p_{ij} = q_{ij} = 0, \\ b_{111}^1 =  \displaystyle
\frac{\frac{1}{2} (x^1)^4 A^1 - (x^1)^2 - 1}{(A^1)^3 B^1}, \;\;
a_{111}^1 = \frac{1}{4} b^1_{111}, \\ b^2_{jkl} = a^2_{jkl} = 0,
\;\; b^1_{jkl} = a^1_{jkl} = 0, \;\; (j, k, l) \neq (1, 1, 1),  \\
f_{11} = g_{11} =  h_{11} =  \frac{1}{4} b^1_{111}, \;\; f_{i2} =
g_{i2} = h_{i2} = 0.
        \end{array}
        \right.
\renewcommand{\arraystretch}{1}
\end{equation}
where
$$
 A^i  = 1 + x^i y^i; \;\; B^i = 1 + \displaystyle
\frac{(x^i)^2}{2}.
$$
It follows from (75), (32),  and (33)
 that {\em the web $(74)$ is isoclinicly geodesic
and not transversally geodesic}.
Thus, {\em this web belongs to the classes {\bf B} and {\bf C}}.
Conditions (75) show that {\em it also belongs to the class ${\bf E_1}$}.
}

\examp{\label{examp:14}
 Consider the web  defined by
\begin{equation}\label{eq:76}
u_3^1 =  x^1 +  y^1 +  x^1 y^2, \;\;\; u_3^2 =  x^2 y^2
\end{equation}
 in a domain of ${\bf R}^4$ where $x^2  \neq 0$ and
$y^2  \neq 0, - 1$.

Using (41)--(46), (26), (28), and (31), for the web (76),
we find that
\begin{equation}\label{eq:77}
\renewcommand{\arraystretch}{1.3}
\left\{
\begin{array}{ll}
\Gamma^1_{12} = -\displaystyle \frac{1}{x^2 (1 + y^2)}, \;\;
\Gamma^1_{jk} = 0, \;\; (j, k) \neq (1, 2), \\
 \Gamma^2_{22} =
-\displaystyle \frac{1}{x^2 y^2}, \;\; \Gamma^2_{jk} = 0, \;\; (j,
k) \neq (2, 2),\\
\omega_1^1 =  -\displaystyle \frac{1}{x^2 (1 +
y^2)} \crazy{\omega}{2}^2, \;\; \omega_1^2 = 0,\\
\omega_2^1 =
-\displaystyle \frac{1}{x^2 (1 + y^2)} \crazy{\omega}{1}^1, \;\;
\omega_2^2 =  -\displaystyle \frac{1}{x^2 y^2}
(\crazy{\omega}{1}^2 + \crazy{\omega}{2}^2),\\
a_1 = 0, \;\; a_2 =
\displaystyle \frac{1}{(1 + y^2) x^2}, \;\;
 q_{22} = \displaystyle
\frac{1}{(1 + y^2)^2 (x^2)^2 y^2},\;\; p_{ij} = q_{1i}  = 0, \\
b^1_{212} = \displaystyle \frac{1}{(1 + y^2)^2 (x^2)^2 y^2}, \;\;
b^1_{jkl} = 0, \;\; (j, k, l) \neq (2, 1, 2);\;\; b^2_{jkl} = 0,\\
a^2_{222} = f_{22} = h_{22} = - \displaystyle \frac{1}{3}
g_{22} = - \displaystyle \frac{1}{4} b^1_{212}, \; \; a^2_{jkl} =
0, \;\; (j, k, l) \neq (2, 2, 2), \\
a^1_{jkl} = 0, \;\; f_{11} =
g_{11} = h_{11} = 0, \;\; f_{12} = g_{12} = h_{12} = 0.
        \end{array}
        \right.
\renewcommand{\arraystretch}{1}
\end{equation}

It follows from (77), (32), and (33)
 that {\em the web $(76)$ is isoclinic and not transversally geodesic}.
Thus, {\em this web belongs to the class {\bf C}}. In addition, by (77),
(20), and (23), {\em it belongs to the classes
{\bf A$_{131}$} and {\bf E$_3$}.}
}

\examp{\label{examp:15}
 Consider the web  defined by
\begin{equation}\label{eq:78}
u_3^1 =  x^1 +  y^1 +  x^1 y^2 + x^2 y^1, \;\;\; u_3^2 =  x^2 y^2
\end{equation}
 in a domain of ${\bf R}^4$ where $x^2  \neq 0, - 1$ and
$y^2  \neq 0, - 1$.

Introduce the following notations:  $ \sigma = \displaystyle \frac{1}{x^2 (1 + y^2)}$ and
$\tau = \displaystyle \frac{1}{y^2 (1 + x^2)}$.
Then
using (41)--(46), (26), (28), and (31), for the web (78),
we find that
\begin{equation}\label{eq:79}
\renewcommand{\arraystretch}{1.3}
\left\{
\begin{array}{ll}
\Gamma^1_{11} = \Gamma^2_{11} = \Gamma^2_{12} = \Gamma^2_{21} = 0,\;\;
\Gamma^1_{12} = -\sigma,\;\; \Gamma^1_{21} = -\tau,\\
\Gamma^2_{22} = -\displaystyle \frac{1}{x^2 y^2},\;\;
\Gamma^1_{22} = \sigma \tau (x^1 + y^1 + x^1 y^2 + x^2 y^1),\\
\omega_1^1 =  -\tau \crazy{\omega}{1}^2 -\sigma \crazy{\omega}{2}^2, \;\;
\omega_1^2 = 0, \\
\omega_2^1 = -\sigma \crazy{\omega}{1}^1 -\tau \crazy{\omega}{2}^1, \;\;
\omega_2^2 = -\displaystyle \frac{1}{x^2 y^2} (\crazy{\omega}{1}^2 +
 \crazy{\omega}{2}^2), \\
a_1 = 0, \;\; a_2 =  \sigma \tau (y^2 -  x^2), \;\; p_{1i} =  q_{1i}  = 0, \;\;
 p_{22} = -\displaystyle\frac{\tau^2}{x^2},\;\;
 q_{22} = \displaystyle
\frac{\sigma^2}{y^2}, \\
b^2_{jkl} = b^1_{111} = b^1_{112} = b^1_{121} = b^1_{211} = b^1_{122} = 0, \\
b^1_{212} = \displaystyle \frac{\sigma^2}{y^2}, \;\;
b^1_{221} =
-\displaystyle \frac{\tau^2}{x^2}, \;\; b^1_{222} = a^1_{222} =
\Gamma^1_{22} a_2, \\
%\displaystyle
%\frac{(y^2 - x^2)(x^1 + y^1 + x^1 y^2
%+ x^2 y^1)}{(x^2)^2 (y^2)^2 (1 + x^2)^2 (1 + y^2)^2},
f_{1i} = g_{1i} = h_{1i} = 0, \;\;
h_{22} = \displaystyle \frac{1}{4} \sigma^2 \tau^2
(y^2 - x^2)(1 - x^2 y^2), \\
f_{22} = -\displaystyle \frac{1}{4} \sigma^2 \tau^2
\bigl[3 x^2 (1 + y^2)^2 + x^2 (1 + x^2)^2\bigr],  \\
g_{22} = \displaystyle \frac{1}{4} \sigma^2 \tau^2
\bigl[3 y^2 (1 + x^2)^2 + y^2 (1 + y^2)^2\bigr].
        \end{array}
        \right.
\renewcommand{\arraystretch}{1}
\end{equation}

It follows from (79), (32), and (33)
 that {\em the web $(78)$ is isoclinic and not transversally geodesic}.
Thus, {\em this web belongs to the class {\bf C}}. In addition, by (79),
(19), and (22), {\em it belongs to the classes
{\bf A$_{131}$} and {\bf E$_4$}.}
}

We will present the results of this section in the following table
in which  we indicate to which classes the webs of our 15  examples
belong.

\vspace*{2mm}
\begin{center}
\begin{tabular}{|c|c|c|c|c|c|c|c|}
  % after \\: \hline or \cline{col1-col2} \cline{col3-col4} ...
  \hline
  Example/Classes & {\bf A} & {\bf B} & {\bf C} & {\bf D} & {\bf E} & {\bf F} & {\bf G} \\ \hline \hline
  1 & {\bf A$_{121}$}   &  & {\bf C$_{11}$} &  & {\bf E$_{71}$} &  {\bf F$_{2}$} &  \\ \hline
  2 & {\bf A$_{121}$}   &  &  &  & {\bf E$_{71}$}  & {\bf F$_{2}$} &  \\ \hline
  3 & {\bf A$_{121}$}   &  & {\bf C$_{1}$} &  &  {\bf E$_{2}$} &  & {\bf G} \\ \hline
  4 & {\bf A$_{1}$}   &  &  & {\bf D$_{231}$} &  {\bf E$_{1}$} & {\bf F$_{1}$} &  \\ \hline
  5 & {\bf A$_{1}$} &  & & {\bf D$_{231}$} &  {\bf E$_{1}$} & {\bf F$_{1}$} &   \\ \hline
  6 & {\bf A$_{1121}$} &  & {\bf C} &  &  {\bf E$_{1}$} &  & {\bf G}  \\ \hline
  7 & {\bf A$_{2}$} &  &  & {\bf D$_{21}$} &  {\bf E$_{8}$} & {\bf F$_{1}$} &  \\ \hline
  %8 & {\bf A$_{2}$} &  & {\bf C} &  &  {\bf E$_{8}$} &  &  \\ \hline
  8 &          & {\bf B} & {\bf C} & &  {\bf E$_{1}$} &  &  \\ \hline
  9 & {\bf A$_{1}$} & {\bf B} &  & {\bf D$_{232}$} &  {\bf E$_{1}$} & {\bf F}$_1$   &  \\ \hline
  10 & {\bf A$_{131}$} & {\bf C} & &  &  {\bf E$_{2}$} &  &  {\bf G}  \\ \hline
  11 & {\bf A$_{131}$} &  {\bf C} &  & &  {\bf E$_{3}$} &  &  \\ \hline
  12 & & {\bf B} & {\bf C} &  &  {\bf E$_{1}$} &  &  \\ \hline
  13 & & {\bf B} & {\bf C} &  &  {\bf E$_{1}$} &  &  \\ \hline
  14 & {\bf A$_{131}$} &  & {\bf C} & &  {\bf E$_{3}$} &  &  \\ \hline
 % 16 & {\bf A$_{1121}$} &  & {\bf C} & &  {\bf E$_{2}$} &  &  \\ \hline
  15 & {\bf A$_{131}$} &  & {\bf C} & &  {\bf E$_{4}$} &  &  \\ \hline \hline
\end{tabular}
\end{center} 
\vspace*{2mm}

Thus the examples considered in this section
 prove  existence of all webs indicated in this table.
Moreover, this proves  existence of more general webs than
those indicated in the table. For example,  existence of
{\bf D$_{231}$} proves  existence of
{\bf D$_{23}$}, {\bf D$_{2}$}, and {\bf D}.

\noindent {\em Author's address}:

\noindent Vladislav V. Goldberg\\ Department of Mathematics\\ New
Jersey Institute of Technology\\ Newark, N.J. 07102, U.S.A.

\vspace*{2mm}

\noindent {\em Author's e-mail address}: vlgold@m.njit.edu

\end{document}